\documentclass[leqno]{article}

\def\diam{\mathop{\rm diam}}
\def\dist{\mathop{\rm dist}}
\def\Lip{\mathop{\rm Lip}}
\def\re{\mathop{\rm Re}}
\def\sign{\mathop{\rm sign}}
\def\span{\mathop{\rm span}}
\def\sym{\mathop{\rm Sym}}
\def\tr{\mathop{\rm tr}}

\newtheorem{theorem}{Theorem}
\newtheorem{lemma}[theorem]{Lemma}
\newtheorem{proposition}[theorem]{Proposition}
\newtheorem{definition}[theorem]{Definition}
\newtheorem{corollary}[theorem]{Corollary}

\newcommand{\begintheorem}{\addtocounter{equation}{1}\begin{theorem}}
\newcommand{\beginlemma}{\addtocounter{equation}{1}\begin{lemma}}
\newcommand{\beginproposition}{\addtocounter{equation}{1}\begin{proposition}}
\newcommand{\begindefinition}{\addtocounter{equation}{1}\begin{definition}}
\newcommand{\begincorollary}{\addtocounter{equation}{1}\begin{corollary}}

\begin{document}

\title{Notes on matrices and calculus}

\author{Stephen William Semmes \\
	Rice University \\
	Houston, Texas}

\date{}

\maketitle

\renewcommand{\thefootnote}{}

\footnotetext{These informal notes are dedicated to Carlos
``Charles'' Tomei, whose influence is clearly visible.  The author
would also like to thank various other people for their comments and
suggestions, and in particular a topics course concerning this
material is taking place at Rice University in the Fall semester of
2003, with the participation of Frank Jones, Amanda Knecht, James
Petersen, and Jun Zhang.}

\tableofcontents

\section{Matrices and linear transformations}
\label{matrices and linear transformations}

	As usual, ${\bf R}$ and ${\bf C}$ denote the real and complex
numbers, respectively.  If $z = x + i \, y$ is a complex number, with
$x, y \in {\bf R}$, then the \emph{complex conjugate} of $z$ is
denoted $\overline{z}$ and defined by
\begin{equation}
	\overline{z} = x - i \, y.
\end{equation}
Notice that
\begin{equation}
	\overline{z + w} = \overline{z} + \overline{w}
\end{equation}
and
\begin{equation}
	\overline{z \, w} = \overline{z} \, \overline{w}
\end{equation}
for complex numbers $z$, $w$.

	If $m$, $n$ are positive integers, we shall denote by
$\mathcal{L}({\bf R}^m, {\bf R}^n)$ the space of real-linear mappings
from ${\bf R}^m$ to ${\bf R}^n$, and by $\mathcal{L}({\bf C}^m, {\bf
C}^n)$ the space of complex-linear mappings from ${\bf C}^m$ to ${\bf
C}^n$.  In the special case where $m = n$, we may simply write
$\mathcal{L}({\bf R}^n)$, $\mathcal{L}({\bf C}^n)$, respectively.
Also when $m = n$, we write $I$ for the identity mapping on ${\bf
R}^n$ or ${\bf C}^n$, as appropriate.

	Using the standard basis for real and complex Euclidean
spaces, linear transformations can be identified with matrices in the
usual manner.  Let us write ${\bf M}_r(m,n)$ and
${\bf M}_c(m,n)$ for the spaces of $m \times n$ real and complex
matrices, respectively.  Thus $\mathcal{L}({\bf R}^m, {\bf R}^n)$,
$\mathcal{L}_c({\bf C}^m, {\bf C}^n)$ can be identified with
${\bf M}r(m,n)$, ${\bf M}_c(m,n)$, respectively, and
in particular addition and scalar multiplication of linear
transformations corresponds to componentwise addition and scalar
multiplication of matrices.

	When $m = n$ we write ${\bf M}_r(n)$ and ${\bf M}_c(n)$ for
the spaces of $n \times n$ real and complex matrices, respectively.
Two elements of ${\bf M}_r(n)$ or of ${\bf M}_c(n)$ can be multiplied
in the customary manner of ``matrix multiplication'', which
corresponds exactly to composition of the associated linear
transformations on ${\bf R}^n$ or ${\bf C}^n$.  The matrix associated
to the identity transformation $I$ has $1$'s along the diagonal and
$0$'s elsewhere, and the product of this matrix with another matrix
gives back that other matrix, just as the composition of the identity
transformation with another transformation gives back that other
transformation.

	If $x = (x_1, \ldots, x_n)$, $y = (y_1, \ldots, y_n)$ are
elements of ${\bf R}^n$, then their \emph{inner product} is denoted
$\langle x, y \rangle$ and is defined by
\begin{equation}
	\langle x, y \rangle = \sum_{j=1}^n x_j \, y_j.
\end{equation}
In the complex case, the inner product of two vectors $z = (z_1,
\ldots, z_n)$ and $w = (w_1, \ldots, w_n)$ is also denoted $\langle z,
w \rangle$ and is defined by
\begin{equation}
	\langle z, w \rangle = \sum_{j=1}^n z_j \, \overline{w_j}.
\end{equation}
In both cases the standard Euclidean norm of an element $v$ of ${\bf
R}^n$ or ${\bf C}^n$ is denoted $|v|$ and is defined to be the
nonnegative real number such that
\begin{equation}
	|v|^2 = \langle v, v \rangle.
\end{equation}

	Given a linear transformation $T$ on ${\bf R}^n$ or ${\bf C}^n$,
there is a unique linear transformation $T^*$ on the same space
such that
\begin{equation}
	\langle T^*(v), w \rangle = \langle v, T(w) \rangle
\end{equation}
for all $v$, $w$ in ${\bf R}^n$ or ${\bf C}^n$, as appropriate.  This
linear transformation $T^*$ is called the \emph{adjoint} of $T$.  In
the real case, the matrix associated to $T^*$ is the \emph{transpose}
of the matrix associated to $T$, which is to say that the $(j, l)$
component of the matrix associated to $T^*$ is equal to the $(l, j)$
component of the matrix associated to $T$, $1 \le j, l \le n$, and in
the complex case the matrix associated to $T^*$ can be obtained by
taking the complex conjugates of the entries of the transpose of the
matrix associated to $T$.

	A linear transformation $T$ on ${\bf R}^n$ or ${\bf C}^n$ is
said to be \emph{self-adjoint} if $T = T^*$, and the space of
self-adjoint linear transformations on ${\bf R}^n$, ${\bf C}^n$ is
denoted $\mathcal{S}({\bf R}^n)$, $\mathcal{S}({\bf C}^n)$,
respectively.  The identity transformation is self-adjoint, and if
$T_1$, $T_2$ are elements of $\mathcal{S}({\bf R}^n)$ or of
$\mathcal{S}({\bf C}^n)$, and if $r_1$, $r_2$ are real numbers, then
the linear combination
\begin{equation}
	r_1 \, T_1 + r_2 \, T_2
\end{equation}
is also an element of $\mathcal{S}({\bf R}^n)$ or $\mathcal{S}({\bf
C}^n)$, respectively.  Note that it is important to use real numbers
as scalars here even if one is working with linear transformations
on ${\bf C}^n$.

	Let us write $\mathcal{S}_r(n)$, $\mathcal{S}_c(n)$ for the
spaces of real and complex $n \times n$ matrices, respectively, that
correspond to self-adjoint linear transformations.  Thus
$\mathcal{S}_r(n)$ consists of the matrices in ${\bf M}_r(n)$
which are \emph{symmetric}, in the sense that the $(j, l)$ and $(l,
j)$ entries are equal to each other.  Similarly, $\mathcal{S}_c(n)$
consists of the matrices in ${\bf M}_c(n)$ such that the $(l, j)$
entry is equal to the complex conjugate of the $(j, l)$ entry, and
in particular so that the diagonal or $(j, j)$ entries are real
numbers.

	If $x, y \in {\bf R}^n$, then
\begin{equation}
	\langle y, x \rangle = \langle x, y \rangle,
\end{equation}
while if $z, w \in {\bf C}^n$, then
\begin{equation}
	\langle w, z \rangle = \overline{\langle z, w \rangle}.
\end{equation}
As a consequence, if $T$ is a self-adjoint linear transformation
on ${\bf C}^n$, and $v$ is an element of ${\bf C}^n$, then
\begin{equation}
	\overline{\langle T(v), v \rangle} = \langle v, T(v) \rangle
		= \langle T(v), v \rangle,
\end{equation}
so that $\langle T(v), v \rangle$ is a real number.  Conversely, if
$T$ is a linear transformation on ${\bf C}^n$ and $\langle T(v), v
\rangle$ is a real number for all $v \in {\bf C}^n$, then $T$ is
self-adjoint.

	In fact, in the complex case a linear transformation $T$ on
${\bf C}^n$ can always be expressed as $S_1 + i \, S_2$, where $S_1$,
$S_2$ are self-adjoint linear transformations on ${\bf C}^n$.  Namely,
one can take $S_1 = (T + T^*)/2$ and $S_2 = (T - T^*)/(2i)$.  It is
easy to see that $\langle T(v), v \rangle$ is real for all $v \in {\bf
C}^n$ if and only if $\langle S_2(v), v \rangle = 0$ for all $v \in
{\bf C}^n$.

	In both the real and complex cases we have the following fact.
Suppose that $S$ is a self-adjoint linear transformation on ${\bf R}^n$
or ${\bf C}^n$ such that
\begin{equation}
	\langle S(v), v \rangle = 0
\end{equation}
for all $v$ in ${\bf R}^n$ or ${\bf C}^n$, as appropriate.  Then $S$
is equal to the zero linear transformation.

	More generally, suppose that $S$ is a self-adjoint linear
transformation on ${\bf R}^n$ or ${\bf C}^n$, and that $v$ is an
element of ${\bf R}^n$ or ${\bf C}^n$ such that $|v| = 1$ and
\begin{equation}
	\langle S(w), w \rangle
\end{equation}
is maximized, or minimized, or has a critical point at $v$, as a
function on the unit sphere, which consists of the vectors $w$ such
that $|w| = 1$.  As in vector calculus, one can check that $v$ is an
\emph{eigenvector} for $S$, in the sense that there is a scalar
$\lambda$ such that
\begin{equation}
	S(v) = \lambda \, v.
\end{equation}
This scalar $\lambda$ is called the \emph{eigenvalue} of $S$
associated to the eigenvector $v$, and for a self-adjoint linear
transformation it is easy to verify that the eigenvalues must be real
numbers, even in the complex case.

	This is the computation used in a standard proof of the fact
that self-adjoint linear operators on ${\bf R}^n$ or ${\bf C}^n$ can
be diagonalized in an orthonormal basis.  In other words, if $S$ is a
self-adjoint linear transformation on ${\bf R}^n$ or ${\bf C}^n$, then
there are eigenvectors $v_1, \ldots, v_n$ for $S$ which are
orthonormal in the sense that
\begin{equation}
	\langle v_j, v_l \rangle = 0
\end{equation}
when $j \ne l$ and
\begin{equation}
	\langle v_j, v_j \rangle = 1
\end{equation}
for each $j$.  Let us also mention that if $S$ is a self-adjoint
linear transformation on ${\bf R}^n$, ${\bf C}^n$ and $v$ is an
eigenvector for $S$, and if $w$ is another vector in ${\bf R}^n$,
${\bf C}^n$ which is orthogonal to $v$ in the sense that
\begin{equation}
	\langle v, w \rangle = 0,
\end{equation}
then $S(w)$ is also orthogonal to $v$.

	A self-adjoint linear transformation $T$ on ${\bf R}^n$ or
${\bf C}^n$ is said to be \emph{nonnegative} if
\begin{equation}
	\langle T(v), v \rangle \ge 0
\end{equation}
for all $v$ in ${\bf R}^n$ or ${\bf C}^n$, as appropriate.  This is
equivalent to the condition that the eigenvalues of $T$ be nonnegative
real numbers.  If $T_1$, $T_2$ are nonnegative self-adjoint linear
transformations on ${\bf R}^n$ or on ${\bf C}^n$ and $r_1$, $r_2$ are
nonnegative real numbers, then
\begin{equation}
	r_1 \, T_1 + r_2 \, T_2
\end{equation}
is also a nonnegative self-adjoint linear transformation.

	A linear transformation $A$ on ${\bf R}^n$ or ${\bf C}^n$ is
said to be \emph{invertible} if there is another linear transformation
$B$ on ${\bf R}^n$ or ${\bf C}^n$, as appropriate, such that
\begin{equation}
	A \circ B = B \circ A = I.
\end{equation}
It is easy to check that if $B$ is a mapping on ${\bf R}^n$ or ${\bf
C}^n$ which is the inverse of $A$ as a mapping, then $B$ must also be
linear, so that $A$ is invertible as a linear mapping.  The inverse of
a linear transformation $A$ is unique when it exists, and is denoted
$A^{-1}$.

	The \emph{kernel} of a linear transformation $A$ on ${\bf
R}^n$ or ${\bf C}^n$ is the set of vectors $v$ in ${\bf R}^n$ or ${\bf
C}^n$, as appropriate, such that $A(v) = 0$. The kernel of a linear
transformation is automatically a linear subspace, which means that it
contains the vector $0$, the sum of two elements of the kernel again
lies in the kernel, and any scalar multiple of a vector in the kernel
is also an element of the kernel.  The kernel of a linear
transformation is said to be \emph{trivial} if it contains only the
vector $0$.

	If a linear transformation is invertible, then its kernel is
trivial.  Conversely, if $A$ is a linear transformation on ${\bf R}^n$
or ${\bf C}^n$ whose kernel is trivial, then $A$ is invertible.  This
is a well known fact from linear algebra, and similarly $A$ is
invertible if and only if it maps ${\bf R}^n$ or ${\bf C}^n$ onto
itself, as appropriate.

	The statement that a linear transformation $A$ on ${\bf R}^n$
or ${\bf C}^n$ is nontrivial is equivalent to the statement that $A$
has a nonzero eigenvector with eigenvalue equal to $0$.  More
generally, a scalar $\lambda$ is an eigenvalue for a linear
transformation $A$ if and only if the linear transformation
\begin{equation}
	A - \lambda \, I
\end{equation}
has a nontrivial kernel.  For the record, a scalar $\lambda$ is
considered to be an eigenvalue of a linear transformation $A$ only
when there is a \emph{nonzero} eigenvector for $A$ with eigenvalue
$\lambda$.

	If $A_1$, $A_2$ are invertible linear transformations on ${\bf
R}^n$ or on ${\bf C}^n$, then the composition $A_1 \circ A_2$ is also
invertible.  In this case we have that
\begin{equation}
	(A_1 \circ A_2)^{-1} = A_2^{-1} \circ A_1^{-1}.
\end{equation}
Conversely, if $A_1$ and $A_2$ are linear transformations on ${\bf
R}^n$ or on ${\bf C}^n$ such that $A_1 \circ A_2$ is invertible, then
$A_1$ and $A_2$ are each invertible themselves, because $A_1$ maps
${\bf R}^n$ or ${\bf C}^n$ onto itself, as appropriate, and $A_2$ has
trivial kernel.

	Suppose that $T_1$, $T_2$ are linear operators on ${\bf R}^n$
or on ${\bf C}^n$.  One can check that
\begin{equation}
	(T_1 \circ T_2)^* = T_2^* \circ T_1^*.
\end{equation}
If $T$ is an invertible linear operator on ${\bf R}^n$ or
${\bf C}^n$, then $T^*$ is also invertible, with
\begin{equation}
	(T^*)^{-1} = (T^{-1})^*,
\end{equation}
and in particular the inverse of an invertible self-adjoint linear
operator is also self-adjoint.

	A self-adjoint linear operator $A$ on ${\bf R}^n$ or ${\bf
C}^n$ is said to be \emph{positive-definite} if
\begin{equation}
	\langle A(v), v \rangle > 0
\end{equation}
for all nonzero vectors $v$.  Thus a positive-definite self-adjoint
linear operator is invertible, because it has trivial kernel, and one
can check that the inverse is also positive-definite.  Also,
a self-adjoint linear transformation is positive-definite if and
only if it is nonnegative and invertible.

	Suppose that $T$ is any linear transformation on ${\bf R}^n$
or ${\bf C}^n$.  Clearly $T^* \circ T$ is self-adjoint, and it is 
nonnegative as well.  Moreover, $T^* \circ T$ is positive definite
if and only if $T$ is invertible.

	If $A$ is a self-adjoint linear transformation on ${\bf R}^n$
or ${\bf C}^n$ which is positive-definite, and if $\alpha$ is a
positive real number, then $\alpha \, A$ is also a self-adjoint
linear transformation which is positive-definite.  If $A_1$, $A_2$
are two self-adjoint linear transformations on ${\bf R}^n$ or on
${\bf C}^n$ which are self-adjoint and nonnegative, and if at least
one of $A_1$, $A_2$ is positive-definite, then the sum $A_1 + A_2$
is a self-adjoint linear transformation which is positive-definite.
In particular, the sum of two self-adjoint linear transformations
which are positive-definite is again positive-definite.

	A linear transformation $T$ on ${\bf R}^n$ or ${\bf C}^n$ is
said to be \emph{orthogonal} or \emph{unitary}, respectively, if $T$
is invertible and
\begin{equation}
	T^{-1} = T^*.
\end{equation}
This is equivalent to saying that
\begin{equation}
	\langle T(v), T(w) \rangle = \langle v, w \rangle
\end{equation}
for all vectors $v$, $w$ in the domain.  In fact, this is equivalent
to
\begin{equation}
	|T(v)| = |v|
\end{equation}
for all vectors $v$ in the domain, as one can show using
\emph{polarization}.

	A linear transformation $A$ on ${\bf R}^n$ or ${\bf C}^n$ is
said to be \emph{anti-self-adjoint} if
\begin{equation}
	A^* = - A.
\end{equation}
Any linear transformation $T$ can be written as $S + A$, with $S$
a self-adjoint linear transformation and $A$ an anti-self-adjoint
linear transformation, simply by taking
\begin{equation}
	S = \frac{T + T^*}{2}, \quad A = \frac{T - T^*}{2}.
\end{equation}
In the complex case a linear transformation is anti-self-adjoint if
and only if it is $i$ times a self-adjoint linear transformation, and
in both the real and complex cases it can be useful to observe that
the square of an anti-self-adjoint operator is self-adjoint, and in
fact it is $-1$ times a nonnegative self-adjoint
operator.

	As above, a subset $L$ of ${\bf R}^n$ or ${\bf C}^n$ is said
to be a \emph{linear subspace} if $0 \in L$, $v, w \in L$ implies $v +
w \in L$, and $v \in L$ implies $\alpha \, v \in L$ for all scalars
$\alpha$, which is to say all real or complex numbers, as appropriate.
The subspace consisting of only the vector $0$ is called the trivial
subspace.  Of course ${\bf R}^n$, ${\bf C}^n$ are linear subspaces of
themselves.

	Suppose that $v_1, \ldots, v_m$ is a finite collection of
vectors in ${\bf R}^n$ or ${\bf C}^n$.  The \emph{span} of $v_1,
\ldots, v_m$ is denoted
\begin{equation}
	\span \{v_1, \ldots, v_m\}
\end{equation}
and is the linear subspace consisting of the vectors of the form
\begin{equation}
	\sum_{j=1}^m \alpha_j \, v_j,
\end{equation}
where $\alpha_1, \ldots, \alpha_m$ are scalars.  A linear subspace $L$
of ${\bf R}^n$ or ${\bf C}^n$ is said to be spanned by a finite
collection of vectors $v_1, \ldots, v_m \in L$ if the span of those
vectors is equal to $L$.

	A finite collection $v_1, \ldots, v_m$ of vectors in ${\bf
R}^n$ or ${\bf C}^n$ is said to be \emph{linearly independent} if a
linear combination $\sum_{j=1}^m \alpha_j \, v_j$ of the $v_j$'s is
equal to the vector $0$ only when the scalars $\alpha_j$ are all equal
to $0$.  This is equivalent to saying that vectors in the span of
$v_1, \ldots, v_m$ are represented in a unique manner as a linear
combination $\sum_{j=1}^m \alpha_j \, v_j$.  A finite collection
$\{v_1, \ldots, v_m\}$ of vectors in ${\bf R}^n$ or ${\bf C}^n$ is
said to be a \emph{basis} for a linear subspace $L$ if $v_1, \ldots,
v_m$ are linearly independent and their span is equal to
$L$.

	A finite collection $v_1, \ldots, v_m$ of vectors in ${\bf
R}^n$ or ${\bf C}^n$ is linearly dependent if there are scalars
$\alpha_1, \ldots, \alpha_m$, with $\alpha_j \ne 0$ for at least one
$j$, such that
\begin{equation}
	\sum_{j=1}^m \alpha_j \, v_j = 0.
\end{equation}
In this case one can reduce the collection to a smaller one with the
same span, at least if we consider the trivial subspace to be the 
span of the empty collection of vectors.  Assuming that at least one
of the vectors is nonzero, we can repeat the process to obtain a
nonempty subcollection of vectors which is linearly independent
and has the same span.

	A basic result from linear algebra states that if $L$ is a
linear subspace of ${\bf R}^n$ or ${\bf C}^n$ which is spanned by a
collection of $m$ vectors, then every linearly independent collection
of vectors in $L$ has less than or equal to $m$ elements.  This comes
down to the fact that a system of $l$ homogeneous linear equations
with more than $l$ variables always has a nontrivial solution.
One can turn this around and say that if $L$ contains a set of $k$
linearly independent vectors, then any collection of vectors which
spans $L$ has at least $k$ elements.

	The \emph{standard basis} for ${\bf R}^n$ or ${\bf C}^n$ is
the collection of $n$ vectors, each of which has exactly one component
equal to $1$ and the others equal to $0$.  It is easy to see that this
is a basis, which is to say that it is linearly independent and spans
the whole space.  Also, every linear subspace of ${\bf R}^n$ or ${\bf
C}^n$ is spanned by a finite collection of vectors, and hence has a
basis, using the empty collection of vectors for the trivial subspace.

	The \emph{dimension} of a linear subspace of ${\bf R}^n$ or
${\bf C}^n$ is equal to the number of elements of a basis in the
subspace.  By the earlier remarks this number is the same for each
basis.  The dimension can also be described as the maximum number of
linearly independent vectors in the subspace, or the minimal number of
vectors needed to span the subspace.

	Let $L$ be a linear subspace of ${\bf R}^n$ or ${\bf C}^n$
with dimension $l$.  A collection of $l$ linearly independent vectors
in $L$ also spans $L$, since otherwise one could add a vector in $L$
not in the span of these vectors to get a collection of $l + 1$
linearly independent vectors in $L$.  Similarly, a collection of $l$
vectors in $L$ which spans $L$ is also linearly independent.

	Suppose that $T$ is a linear operator on ${\bf R}^n$ or ${\bf
C}^n$, and that $L$ is a linear subspace of the same space.  In this
event $T(L)$, the image of $L$ under $T$, is also a linear subspace.
If $T$ is invertible, then the dimension of $T(L)$ is equal to the
dimension of $L$.

	A collection of vectors $v_1, \ldots, v_m$ in ${\bf R}^n$ or
${\bf C}^n$ is said to be \emph{orthonormal} if, as before,
\begin{equation}
	\langle v_j, v_k \rangle = 0
\end{equation}
when $j \ne k$ and
\begin{equation}
	|v_j| = 1
\end{equation}
for each $j$.  If $v_1, \ldots, v_m$ is an orthonormal collection of
vectors in ${\bf R}^n$ or ${\bf C}^n$ and $w$ is in their span, so
that
\begin{equation}
	w = \sum_{j=1}^m \alpha_j \, v_j
\end{equation}
for some scalars $\alpha_j$, then
\begin{equation}
	\alpha_j = \langle w, v_j \rangle
\end{equation}
for each $j$, and in particular $v_1, \ldots, v_m$ are linearly
independent.  Also, we have that
\begin{equation}
	|w|^2 = \sum_{j=1}^m |\langle w, v_j \rangle|^2
\end{equation}
in this case.

	Let us recall the \emph{Cauchy--Schwarz inequality}, which
states that if $v$, $w$ are elements of ${\bf R}^n$ or of ${\bf C}^n$,
then
\begin{equation}
	|\langle v, w \rangle| \le |v| \, |w|.
\end{equation}
This can be shown using the fact that
\begin{equation}
	\langle v + \alpha \, w, v + \alpha \, w \rangle 
		= |v + \alpha \, w|^2 \ge 0
\end{equation}
for all scalars $\alpha$.  Using this inequality, one can also show
that
\begin{equation}
	|v + w| \le |v| + |w|,
\end{equation}
which is to say the triangle inequality.

	As before, if $v$, $w$ are two vectors in ${\bf R}^n$ or ${\bf
C}^n$, then we say that $v$, $w$ are \emph{orthogonal} if
\begin{equation}
	\langle v, w \rangle = 0,
\end{equation}
and in this case we write $v \perp w$.  If $v$, $w$ are orthogonal
vectors in ${\bf R}^n$ or in ${\bf C}^n$, then
\begin{equation}
	|v + w|^2 = |v|^2 + |w|^2,
\end{equation}
and $\alpha \, v \perp \beta \, w$ for all scalars $\alpha$, $\beta$.
Conversely, notice that if $v$, $w$ are two vectors in ${\bf R}^n$
such that $|v + w|^2 = |v|^2 + |w|^2$, then $v \perp w$, and if $v$,
$w$ are two vectors in ${\bf C}^n$ such that $|v + \alpha \, w|^2 =
|v|^2 + |w|^2$ for all complex numbers $\alpha$ with $|\alpha| = 1$,
then $v \perp w$.

	Suppose again that $v_1, \ldots, v_m$ is an orthonormal
collection of vectors in ${\bf R}^n$ or ${\bf C}^n$.  If $u$
is any vector in ${\bf R}^n$ or ${\bf C}^n$, then
\begin{equation}
	u' = \sum_{j=1}^m \langle u, v_j \rangle \, v_j
\end{equation}
lies in the span of $v_1, \ldots, v_m$, and one can check that
$u - u'$ is orthogonal to every vector in the linear span of
$v_1, \ldots, v_m$.  In particular,
\begin{equation}
	\langle u - u', u' \rangle = 0,
\end{equation}
and
\begin{equation}
	|u|^2 = |u - u'|^2 + |u'|^2.
\end{equation}

	If $T$ is a linear transformation on ${\bf R}^n$ or ${\bf
C}^n$, then the \emph{trace} of $T$ is denoted $\tr T$ and is defined
to be the sum of the diagonal terms in the standard matrix associated
to $T$.  To be more explicit, let $e_1, \ldots, e_n$ denote the
standard basis for ${\bf R}^n$ or ${\bf C}^n$, as appropriate, so that
$e_j$ has $j$th component equal to $1$ and all other components equal
to $0$.  The trace of a linear transformation $T$ can then be
expressed as
\begin{equation}
	\tr T = \sum_{j=1}^n \langle T(e_j), e_j \rangle.
\end{equation}

	The trace is clearly linear in $T$, so that if $T_1$, $T_2$
are linear transformations on ${\bf R}^n$ or on ${\bf C}^n$ and
$\alpha_1$, $\alpha_2$ are scalars, then
\begin{equation}
	\tr (\alpha_1 \, T_1 + \alpha_2 \, T_2) 
		= \alpha_1 \, \tr T_1 + \alpha_2 \, \tr T_2.
\end{equation}
Another fundamental property of the trace is that
\begin{equation}
	\tr (T_1 \circ T_2) = \tr (T_2 \circ T_1)
\end{equation}
for all linear transformations $T_1$, $T_2$.  This can be verified
in a straightforward manner.

	If $T$ is a linear transformation on ${\bf R}^n$, then
\begin{equation}
	\tr T^* = \tr T,
\end{equation}
while if $T$ is a linear transformation on ${\bf C}^n$, then
\begin{equation}
	\tr T^* = \overline{\tr T}.
\end{equation}
If $A$, $B$ are linear transformations on ${\bf R}^n$ or ${\bf C}^n$,
then
\begin{equation}
	\tr (B^* \circ A) 
	   	= \sum_{j=1}^n \langle B^*(A(e_j)), e_j \rangle	
	   	= \sum_{j=1}^n \langle A(e_j), B(e_j) \rangle.	
\end{equation}
where as usual $e_1, \ldots, e_n$ denotes the standard basis for ${\bf
R}^n$ or ${\bf C}^n$.  This is the same as
\begin{equation}
	\sum_{j=1}^n \sum_{l=1}^n
   \langle A(e_j), e_l \rangle \langle B(e_j), e_l \rangle
\end{equation}
in the real case and
\begin{equation}
	= \sum_{j=1}^n \sum_{l=1}^n
   \langle A(e_j), e_l \rangle \overline{\langle B(e_j), e_l \rangle}
\end{equation}
in the complex case, which is to say that one takes the standard
matrices of $A$, $B$, views them as elements of ${\bf R}^{n^2}$ or
${\bf C}^{n^2}$, as appropriate, and then takes the usual inner
product.

	In particular, if $T$ is a linear transformation on ${\bf
R}^n$ or ${\bf C}^n$, let $\|T\|_{HS}$ be the nonnegative real number
defined by
\begin{equation}
	\|T\|_{HS}^2 = \tr (T^* \circ T) 
		= \sum_{j=1}^n \sum_{k=1}^n |\langle T(e_j), e_k \rangle|^2.
\end{equation}
In other words, $\|T\|_{HS}$ is the same as the usual Euclidean norm
of the standard matrix associated to $T$, and it is also known as the
\emph{Hilbert--Schmit norm} of $T$.  Observe that $\|T\|_{HS} = 0$ if
and only if $T = 0$, $\|\alpha \, T \|_{HS} = |\alpha| \, \|T\|_{HS}$
for all scalars $\alpha$ and all linear transformations $T$, $\|T_1 +
T_2\|_{HS} \le \|T_1\|_{HS} + \|T_2\|_{HS}$ for all linear
transformations $T_1$, $T_2$, and that
\begin{equation}
	|\tr (T_2^* \circ T_1)| \le \|T_1\|_{HS} \, \|T_2\|_{HS}
\end{equation}
for all linear transformations $T_1$, $T_2$.

	If $T$ is a linear transformation on ${\bf R}^n$ or ${\bf
C}^n$, then the \emph{operator norm} of $T$ is denoted $\|T\|_{op}$ and
defined to be the maximum of $|T(v)|$ over all vectors $v$ in the
domain with $|v| = 1$, which exists by the extreme value theorem in
calculus.  In other words,
\begin{equation}
	|T(w)| \le \|T\|_{op} \, |w|
\end{equation}
for all vectors $w$, and $\|T\|_{op}$ is the smallest nonnegative real
number with this property.  One can check that $\|T\|_{op} = 0$ if and
only if $T = 0$, $\|\alpha \, T\|_{op} = |\alpha| \, \|T\|_{op}$ for
all scalars $\alpha$ and all linear transformations $T$, $\|T_1 +
T_2\|_{op} \le \|T_1\|_{op} + \|T_2\|_{op}$ for all linear
transformations $T_1$, $T_2$, and
\begin{equation}
	\|T_1 \circ T_2\|_{op} \le \|T_1\|_{op} \, \|T_2\|_{op}
\end{equation}
for all linear transformations $T_1$, $T_2$.

	Alternatively, $\|T\|_{op}$ can be described as the maximum
of $|\langle T(v), w \rangle|$ over all vectors $v$, $w$ in the
domain such that $|v| = |w| = 1$, which is the same as saying that
\begin{equation}
	|\langle T(v), w \rangle| \le \|T\|_{op} \, |v| \, |w|
\end{equation}
for all vectors $v$, $w$ in the domain, and that $\|T\|_{op}$ is the
smallest nonnegative real number with this property.  In particular,
it follows that
\begin{equation}
	\|T^*\|_{op} = \|T\|_{op}
\end{equation}
for all linear transformations $T$.  It is easy to check as well
that
\begin{equation}
	\|T^*\|_{HS} = \|T\|_{HS}
\end{equation}
for all linear transformations $T$.

	Notice that
\begin{equation}
	|\langle T(e_j), e_k \rangle| \le \|T\|_{op},
\end{equation}
so that the operator norm of $T$ is greater than or equal to the
absolute values of the entries of the standard matrix associated to
$T$.  One can express $\|T\|_{HS}$ by
\begin{equation}
	\|T\|_{HS}^2 = \sum_{j=1}^n |T(e_j)|^2,
\end{equation}
from which it follows that $\|T\|_{HS} \le n^{1/2} \, \|T\|_{op}$.
From this formula it also follows that $\|A \circ B\|_{HS} \le
\|A\|_{op} \, \|B\|_{HS}$, and similarly one has $\|A \circ B\|_{HS}
\le \|A\|_{HS} \, \|B\|_{op}$ for all linear transformations $A$, $B$.

	Suppose that $v_1, \ldots, v_m$ is an orthonormal collection
of vectors in ${\bf R}^n$ or ${\bf C}^m$, and let $L$ denote the span
of this collection.  As we have seen, if $u$ is any vector in ${\bf
R}^n$ or ${\bf C}^n$, as appropriate, then there is a vector $u' \in
L$ such that $u - u'$ is orthogonal to every element of $L$.
These two properties characterize $u'$, since if $u'' \in L$ and
$u - u''$ is orthogonal to every element of $L$, then $u' - u''$
is an element of $L$ and is orthogonal to every element of $L$,
including itself, so that $u' - u'' = 0$.

	In this situation let us write $P_L$ for the linear
transformation on ${\bf R}^n$ or ${\bf C}^n$, as appropriate, which
sends $u$ to $u'$.  This is called the \emph{orhogonal projection} of
${\bf R}^n$ or ${\bf C}^n$, as appropriate, onto $L$.  It is uniquely
determined by $L$, which is to say that it does not depend on the
choice of orthonormal basis for $L$.

	Using these orthogonal projections, one can show that every
orthonormal set of vectors in ${\bf R}^n$ or ${\bf C}^n$ can be
extended to an orthonormal basis, and that every linear subspace of
${\bf R}^n$ or ${\bf C}^n$ has an orthonormal basis.  This is
basically the same as the Gram--Schmit process, in which a collection
of vectors is orthonormalized one step at a time.  In particular,
for every linear subspace $L$ of ${\bf R}^n$ or ${\bf C}^n$ there
is a corresponding orthogonal projection $P_L$, which one can also
check is self-adjoint.

	Let $v_1, \ldots, v_n$ and $w_1, \ldots, w_n$ be orthonormal
bases of ${\bf R}^n$ or of ${\bf C}^n$, respectively.  If $T$ is a
linear transformation on ${\bf R}^n$ or on ${\bf C}^n$, as appropriate,
then one can check that
\begin{equation}
	\|T\|_{HS}^2 = \sum_{j=1}^n |T(v_j)|^2
\end{equation}
and that
\begin{equation}
	\|T\|_{HS}^2 = 
	   \sum_{j=1}^n \sum_{k=1}^n |\langle T(v_j), w_k \rangle|^2.
\end{equation}
In particular, it follows that $\|T\|_{op} \le \|T\|_{HS}$.

	Suppose that $L$ is a nontrivial linear subspace of ${\bf
R}^n$ or ${\bf C}^n$, and that $P_L$ is the corresponding orthogonal
projection onto $L$.  For each vector $u$ in the domain, we have that
$P_L(u)$ and $u - P_L(u)$ are orthogonal to each other, so that
\begin{equation}
	|P_L(u)|^2 + |u - P_L(u)|^2 = |u|^2,
\end{equation}
and one can check that $\|P_L\|_{op} = 1$.  From the remarks in the
previous paragraphs it follows that $\|P_L\|_{HS}$ is equal to the
square root of the dimension of $L$.

	In general, a \emph{projection} on ${\bf R}^n$ or ${\bf C}^n$
is a linear operator $P$ which is an ``idempotent'', which means that
\begin{equation}
	P^2 = P.
\end{equation}
Thus for instance the identity and the operator $0$ are projections,
and in general if $P$ is a projection and $L$ is the image of $P$, so
that $L$ consists of the vectors of the form $P(v)$ for vectors $v$ in
the domain, then $L$ is exactly the set of vectors $w$ such that $P(w)
= w$.  If $P$ is a projection and $v$ is any vector in the domain,
then $P(v)$ lies in the image of $P$ and $v - P(v)$ lies in the kernel
of $P$.

	If $L$ is a linear subspace of ${\bf R}^n$ or ${\bf C}^n$,
then the \emph{orthogonal complement} of $L$ is denoted $L^\perp$ and
defined to be the linear subspace of vectors $v$ such that $v$ is
orthogonal to $w$ for all $w \in L$.  From the earlier remarks it
follows that every vector $u$ in ${\bf R}^n$ or ${\bf C}^n$, as
appropriate, can be written in a unique way as the sum of vectors in
$L$ and $L^\perp$.  One can also check that $(L^\perp)^\perp = L$.

	A projection $P$ on ${\bf R}^n$ or ${\bf C}^n$ with image
$L$ is equal to the orthogonal projection onto $L$ if and only
if the kernel of $P$ is equal to $L^\perp$.  Also, a projection $P$
is an orthogonal projection if and only if $P$ is self-adjoint.
The operator norm of a nonzero projection is automatically greater
than or equal to $1$, and one can check that it is equal to $1$
if and only if the projection is an orthogonal projection.

	Now let us briefly review some aspects of \emph{determinants}.
We begin with some facts about \emph{permutations}.  Fix a positive
integer $n$, and let $\sym (n)$ denote the \emph{symmetric group on
$\{1, \ldots, n\}$} consisting of the permutations on the set $\{1,
\ldots, n\}$ of positive integers from $1$ to $n$, which is to say the
one-to-one mappings from this set onto itself, with composition
mappings as the group operation, and inverses of mappings as inverses
in the group.

	A \emph{transposition} is a permutation $\tau$ on $\{1,
\ldots, n\}$ which interchanges two elements of the set and leaves the
others fixed.  A basic fact is that every element of the symmetric
group can be expressed as a composition of finitely many
transpositions.  Of course such a product is not unique, and another
important result is that the parity of the number of transpositions
used is unique, i.e., it depends only on the original permutation.

	In effect this is the same as saying that the identity
permutation, which fixes all elements of the set, can be expressed as
a composition of an even number of transpositions, and not an odd
number of transpositions.  An element of the symmetric group is said
to be even or odd according to whether it can be expressed as the
composition of an even or odd number of transpositions.  The
composition of two even permutations is even, the composition of two
odd permutations is an odd permutation, the composition of an even and
an odd permutation is an odd permutation, and the inverse of a
permutation $\pi$ has the same type as $\pi$ does.

	Now let $A$ be a linear transformation on ${\bf R}^n$ or
${\bf C}^n$, and let $(a_{j,l})$ denote the corresponding $n \times n$
matrix of real or complex numbers.  The determinant of $A$ is denoted
$\det A$ and is the real or complex number, respectively, given by
\begin{equation}
	\det A = 
		  \sum_{\pi \in \sym (n)} \sign (\pi) \,
			a_{1, \pi(1)} \, a_{2, \pi(2)} \cdots a_{n, \pi(n)},
\end{equation}
where $\sign (\pi)$ is equal to $+1$ or $-1$ according to whether the
permutation $\pi$ is even or odd.  Thus the determinant of $A$ is a
homogeneous polynomial of degree $n$ as a function of the entries
of the matrix $(a_{j, l})$.

	When $n = 1$, the matrix associated to $A$ is really just a 
single number, and the determinant of $A$ is that number.  In general
we have that $\det I = 1$, $\det A^* = \det A$ for all $A$, and
\begin{equation}
	\det (A \circ B) = (\det A) (\det B)
\end{equation}
for all linear transformations $A$, $B$ on ${\bf R}^n$ or on
${\bf C}^n$.  It follows from this that if $A$ is an invertible
linear transformation, then $\det A \ne 0$ and indeed
\begin{equation}
	(\det A)^{-1} = \det (A^{-1}),
\end{equation}
and conversely there is the well-known \emph{Cramer's rule}, which
states that a linear transformation with nonzero determinant is
invertible, with a formula for the inverse in terms of the determinant
of the linear transformation and determinants of submatrices of the
associated matrix.

	Let $v_1, \ldots, v_n$ be a basis for ${\bf R}^n$ or ${\bf
C}^n$, and let $A$ be a linear transformation on ${\bf R}^n$ or ${\bf
C}^n$.  It is easy to see that $A$ is uniquely determined by its
values on $v_1, \ldots, v_n$, and conversely that if $w_1, \ldots,
w_n$ is any other collection of $n$ vectors in ${\bf R}^n$ or ${\bf
C}^n$, as appropriate, then there is a linear transformation $A$ such
that $A(v_j) = w_j$ for each $j$.  Also, $A$ is an invertible linear
transformation if and only if $A(v_1), \ldots, A(v_n)$ is a basis too.

	For any choice of basis for ${\bf R}^n$ or ${\bf C}^n$, there
is a natural correspondence between linear transformations on ${\bf
R}^n$ or ${\bf C}^n$ and matrices with real or complex entries,
respectively, in such a way that the diagonal matrices correspond
exactly to linear transformations for which the vectors in the basis
are eigenvectors.  Of course for any two choices of bases there is
an invertible linear transformation which sends one basis to the other.
For a single linear transformation, one gets two matrices associated
to the two bases, and these two matrices are related by conjugation.

	In particular, for a linear transformation $A$ on ${\bf R}^n$
or on ${\bf C}^n$ and a choice of basis $v_1, \ldots, v_n$, one gets
a matrix associated to this linear transformation and basis, and one
can take the trace or determinant of this matrix.  It follows from the
basic identities for the trace and determinant that the trace of this
matrix is the same as for the matrix associated to any other choice of
basis.  As a special case, if the linear transformation is diagonalizable,
in the sense that there is a basis of eigenvectors, then the trace is
the same as the sum of the corresponding $n$ eigenvalues, and the
determinant is equal to the product of the eigenvalues.

	As another basic example, if $P$ is a projection on ${\bf
R}^n$ or on ${\bf C}^n$ whose image is a linear subspace $L$ of
dimension $l$, then the trace of $P$ is equal to $l$.

	Now let us look at exponentials, beginning with exponentiation
of real numbers.  The exponential function is denoted $\exp(x)$ and
can be defined by the series expansion
\begin{equation}
	\exp(x) = \sum_{n = 0} \frac{x^n}{n!}.
\end{equation}
Here $x^n$ is interpreted as being equal to $1$ when $n = 0$, even if
$x = 0$, and $n!$ is ``$n$ factorial'', the product of the positive
integers from $1$ to $n$, which is also interpreted as being equal to
$1$ when $n = 0$.

	By standard results, this series converges for all $x \in {\bf
R}$, and converges absolutely, and it also converges uniformly on
bounded subsets of ${\bf R}$.  The sum defines a real-valued function
on the real line which is continuous and has continuous derivatives of
all orders, with the derivatives being given by the series obtained by
differentiating this one term by term.  In this case we have the
well-known identity
\begin{equation}
	\exp'(x) = \exp(x),
\end{equation}
i.e., the derivative of the exponential function is itself.

	A related identity is
\begin{equation}
	\exp(x + y) = \exp(x) \, \exp(y).
\end{equation}
Formally this can be derived by multiplying the series for $\exp(x)$
and $\exp(y)$, group terms of total degree $n$, and using the binomial
theorem to identify them with the terms of $\exp(x + y)$.  Convergence
issues can be handled using absolute convergence of the series
involved, by standard arguments.

	Clearly $\exp(x) \ge 1$ when $x \ge 0$.  From the multiplicative
identity it follows that $\exp(x) \ne 0$ for all $x \in {\bf R}$,
and in fact that
\begin{equation}
	\exp(-x) = \frac{1}{\exp(x)}.
\end{equation}
It follows that $\exp(x) > 0$ for all $x \in {\bf R}$, and hence the
derivatives of $\exp(x)$ are all positive as well, so that $\exp(x)$
is strictly increasing and strictly convex in particular.

	Next we consider complex numbers.  That is, we define
$\exp(z)$ for $z \in {\bf C}$ by the same series as before, namely,
\begin{equation}
	\exp(z) = \sum_{n=0}^\infty \frac{z^n}{n!}.
\end{equation}
This series converges absolutely for all $z \in {\bf C}$,
it converges uniformly on bounded subsets of ${\bf C}$, and
it is continuously differentiable of all orders.

	Again we have the identities
\begin{equation}
	\exp'(z) = \exp(z)
\end{equation}
and
\begin{equation}
	\exp(z + w) = \exp(z) \, \exp(w)
\end{equation}
for all $z, w \in {\bf C}$.  The meaning of the differential equation
for $\exp(z)$ is that $\exp(z)$ is a holomorphic function of $z$ whose
complex derivative is equal to $\exp(z)$.  To put it another way, the
differential of $\exp(z)$ at a point $z$ is given by multiplication by
$\exp(z)$, so that
\begin{equation}
	\exp(z + h) = \exp(z) + \exp(z) \cdot h + O(h^2).
\end{equation}

	If $z = x + i \, y$, with $x, y \in {\bf R}$, then
\begin{equation}
	\exp(z) = \exp(x) \, (\cos(y) + i \, \sin(y)).
\end{equation}
This is a well-known and striking formula, which can be seen by
writing out the series expansions for the real and imaginary parts of
$\exp(i y)$ and comparing them with the usual series expansions for
the cosine and sine.  Also, as a complex-valued function of a real
variable, we have that
\begin{equation}
	\frac{d}{dy} \exp(i y) = i \exp(i y)
\end{equation}
and hence
\begin{equation}
	\frac{d^2}{dy^2} \exp(i y) = - \exp(i y),
\end{equation}
which correspond to standard formulas for the derivatives of the
cosine and sine, including the second-order differential equations
that they satisfy.

	It is clear from the series expansion that
\begin{equation}
	\overline{\exp(z)} = \exp(\overline{z})
\end{equation}
for all $z \in {\bf C}$.  In particular, if $z = x + i \, y$ with
$x, y \in {\bf R}$, then
\begin{equation}
	|\exp(z)| = \exp(x).
\end{equation}
The special case
\begin{equation}
	|\exp(i y)| = 1
\end{equation}
for all $y \in {\bf R}$ corresponds to the usual identity 
$\cos(y)^2 + \sin(y)^2 = 1$.

	Fix a positive integer $n$, and suppose that $A$ is a linear
transformation on ${\bf R}^n$ or on ${\bf C}^n$.  We would like to
define $\exp(A)$ by the series
\begin{equation}
	\exp(A) = \sum_{k=0}^\infty \frac{A^k}{k!},
\end{equation}
where now $A^k$ denotes the $k$-fold composition of $A$ as a linear
transformation, interpreted as being the identity operator $I$ when $k
= 0$.  The convergence of this series can be defined in terms of the
convergence of the entries of the corresponding matrices, and as
before we have absolute convergence for all linear transformations
$A$, uniform convergence on bounded sets of such linear
transformations, and that the exponential defines a continuous
function from linear transformations to themselves which is
continuously differentiable of all orders.

	A convenient way to look at absolute convergence of series of
linear transformations is in terms of convergence of the corresponding
series of operator norms.  In this case we have such convergence, 
because
\begin{equation}
	\sum_{k=0}^\infty \frac{\|A^k\|_{op}}{k!}
		\le \sum_{k=0}^\infty \frac{\|A\|_{op}^k}{k!}.
\end{equation}
In particular, note that
\begin{equation}
	\|\exp(A)\|_{op} \le \exp(\|A\|_{op}),
\end{equation}
where the right side refers to the exponential of the operator
norm of $A$ as a real number.

	If $A$ and $B$ are linear transformations on ${\bf R}^n$
or on ${\bf C}^n$ which commute, then we still have that
\begin{equation}
	\exp(A + B) = \exp(A) \circ \exp(B),
\end{equation}
for essentially the same reasons as before.  Of course 
\begin{equation}
	\exp (0) = I,
\end{equation}
and for any linear transformation $A$ we have that $\exp(A)$ is
invertible, with
\begin{equation}
	(\exp(A))^{-1} = \exp(-A)
\end{equation}
and thus
\begin{equation}
	\|(\exp(A))^{-1}\|_{op} \le \exp \|A\|_{op}.
\end{equation}
Furthermore,
\begin{equation}
	(\exp(A))^* = \exp(A^*),
\end{equation}
which is to say that the adjoint of the exponential of $A$
is equal to the exponential of the adjoint of $A$.

	If $A$ is a linear transformation on ${\bf R}^n$ or on ${\bf
C}^n$, and if $T$ is another linear transformation on ${\bf R}^n$ or
on ${\bf C}^n$, as appropriate, which commutes with $A$, then $T$ also
commutes with $\exp(A)$, and the directional derivative of $\exp(A)$
at $A$ in the direction of $T$ is given by multiplication by
$\exp(A)$, so that
\begin{equation}
	\exp(A + T) = \exp(A) + \exp(A) \, T + O(\|T\|_{op}^2).
\end{equation}
In particular, if $A$ is a linear transformation on ${\bf R}^n$ or on
${\bf C}^n$ and we put
\begin{equation}
	E_A(t) = \exp(t \, A),
\end{equation}
viewed as a function on the real line with values in linear transformations
on ${\bf R}^n$ or on ${\bf C}^n$, as appropriate, then this function is
continuously differentiable of all orders and satisfies
\begin{equation}
	\frac{d}{dt} E_A(t) = A \circ E_A(t),
		\quad E_A(0) = I.
\end{equation}
These conditions characterize $E_A(t)$ uniquely, by standard results
about ordinary differential equations.

	What about the determinant of the exponential of a linear
transformation?  Notice first that the differential of the determinant
as a function on linear transformations on ${\bf R}^n$ or on ${\bf
C}^n$ and evaluated at the identity transformation is given by the
trace.  That is, if $T$ is any linear transformation on ${\bf R}^n$
or on ${\bf C}^n$, then
\begin{equation}
	\det (I + T) = 1 + \tr T + O(\|T\|_{op}^2),
\end{equation}
and of course this is just a simple algebraic statement, since the
determinant of a linear transformation $A$ is a polynomial in the
entries of the matrix associated to $A$.

	This implies that
\begin{equation}
	\frac{d}{dt} \det (\exp (t \, A)) = (\tr A) \det (\exp(t \, A)).
\end{equation}
More precisely, at $t = 0$ this follows exactly from the remarks of
the preceding paragraph.  In general, for any two real numbers $r$, $s$,
we have that
\begin{equation}
	\exp ((r + s) \, A) = \exp (r \, A) \circ \exp(s \, A),
\end{equation}
and this permits one to derive the formula for the derivative at any
real number $t$ from the special case of $t = 0$.

	Of course the determinant of $\exp (t \, A)$ at $t = 0$
is equal to $1$, and it follows that
\begin{equation}
	\det (\exp(t \, A)) = \exp (t \, \tr A).
\end{equation}
The trace of $A$ is a real or complex number, and the right side is
the usual exponential of a scalar.  We may as well apply this to $t =
1$ and say that
\begin{equation}
	\det (\exp(A)) = \exp (\tr A).
\end{equation}

	Suppose that $A$ is a linear transformation on ${\bf R}^n$ or
on ${\bf C}^n$ and that $v$ is a vector in ${\bf R}^n$ or in ${\bf
C}^n$, as appropriate.  Set
\begin{equation}
	h(t) = \exp (t \, A)(v)
\end{equation}
viewed as a function from the real line into ${\bf R}^n$ or ${\bf
C}^n$, i.e., where for each $t$ we let $h(t)$ be the image of $v$
under $\exp(t \, A)$.  As before we have that $h'(t) = A(h(t))$ and
that $h(0) = 0$, and $h(t)$ is characterized by these properties by
standard results about ordinary differential equations.

	Assume further that $v$ is an eigenvector for $A$ with eigenvalue
$\lambda$, so that
\begin{equation}
	A(v) = \lambda \, v,
\end{equation}
where $\lambda$ is a scalar.  In this case
\begin{equation}
	\exp (t \, A)(v) = \exp(t \, \lambda) \, v,
\end{equation}
where $\exp (t \, \lambda)$ is the usual exponential mapping for
scalars.  One can see this either from the series expansion for the
exponential or from the characterization in terms of ordinary
differential equations.

	It may be that $A$ is diagonalizable, so that there is a basis
of eigenvectors for $A$.  The elements of this basis are then
eigenvectors for $\exp (t \, A)$ too, with the eigenvalues for the
exponential being given by the exponentials of the corresponding
eigenvalues, as in the previous paragraph.  In other words, the
exponential is then also diagonalizable, and by the same basis as for
$A$ itself.

	For that matter, suppose that $A$ is a linear transformation
on ${\bf R}^n$ or on ${\bf C}^n$, and that $L$ is a linear subspace 
of the same space which is invariant under $A$.  This means that
\begin{equation}
	A(L) \subseteq L,
\end{equation}
which is to say that $A(v) \in L$ for all $v \in L$.  In this event
$L$ is invariant under $\exp(t \, A)$ for all $t$ as well, as one can
see from either the series expansion for the exponential or the
characterization in terms of ordinary differential equations.

	Next we review some aspects of spectral theory of matrices.
If $A$ is a linear transformation on ${\bf C}^n$, then the
\emph{characteristic polynomial} associated to $A$ is defined by
\begin{equation}
	q_A(z) = \det (z \, I - A).
\end{equation}
Thus $q_A(z)$ is a polynomial of degree $n$ whose leading coefficient
is equal to $1$ and which vanishes exactly at the eigenvalues of $A$.

	The fundamental theorem of algebra states that every
nonconstant polynomial on the complex numbers has a root.  As a
result, every linear transformation on ${\bf C}^n$ has at least one
eigenvalue.  Recall as well that every nonconstant polynomial on the
complex numbers can be factored as a nonzero complex number times a
product of linear factors of the form $(z - \alpha)$, $\alpha \in {\bf
C}$.

	If $p(z)$ is a polynomial,
\begin{equation}
	p(z) = c_m \, z^m + c_{m-1} \, z^{m-1} + \cdots + c_0,
\end{equation}
$c_0, \ldots, c_m \in {\bf C}$, and $A$ is a linear transformation on
${\bf C}^n$, then we can define $p(A)$ to be the linear transformation
on ${\bf C}^n$ given by
\begin{equation}
	p(A) = c_m \, A^m + c_{m-1} \, A^{m-1} + \cdots + c_0 \, I.
\end{equation}
Notice that if $p_1$, $p_2$ are polynomials, so that the sum
$p_1 + p_2$ and the product $p_1 \, p_2$ are also polynomials, then
we have that
\begin{equation}
	(p_1 + p_2)(A) = p_1(A) + p_2(A)
\end{equation}
and
\begin{equation}
	(p_1 \, p_2)(A) = p_1(A) \, p_2(A) = p_2(A) \, p_1(A).
\end{equation}
Moreover, the composition $p_1 \circ p_2$ is also a polynomial,
and $(p_1 \circ p_2)(A) = p_1(p_2(A))$.

	If $A$ is a linear transformation on ${\bf C}^n$, $v$ is a
vector in ${\bf C}^n$ which is an eigenvector for $A$ with eigenvalue
$\lambda$, and $p(z)$ is a polynomial, then $v$ is also an eigenvector
for $p(A)$, with eigenvalue $p(\lambda)$.  Conversely, if $A$ is a
linear transformation on ${\bf C}^n$ and $h(z)$ is a polynomial such that
$h(\lambda) \ne 0$ for all eigenvalues $\lambda$ of $A$, then $h(A)$
is invertible.  As a consequence, for a linear transformation $A$
on ${\bf C}^n$ and a polynomial $p(z)$, if a complex number $\mu$
is an eigenvalue of $p(A)$, then there is an eigenvalue $\lambda$
of $A$ such that $p(\lambda) = \mu$.

	The famous Cayley--Hamilton theorem states that for a
linear transformation $A$ on ${\bf C}^n$ and its characteristic
polynomial $q_A(z)$ as above, we have that
\begin{equation}
	q_A(A) = 0.
\end{equation}
It follows that
\begin{equation}
	p(A) = 0
\end{equation}
whenever $p(z)$ is a polynomial which can be expressed as the product
of the characteristic polynomial $q_A(z)$ and another polynomial.
This holds when $p(z)$ vanishes at each eigenvalue of $A$, and to at
least the same order as $q_A$ does.

	In particular, we have that $A^n$ can be expressed as a linear
combination of $A^k$, $1 \le k \le n-1$, and the identity operator.
By repeating this, every positive integer power of $A$ can be expressed
as a linear combination of $A^k$, $1 \le k \le n - 1$, and the
identity operator.  To put it another way, for each polynomial
$p(z)$ there is a polynomial $\widetilde{p}(z)$ of degree at most
$n - 1$ such that $p(A) = \widetilde{p}(A)$.

	Also, the exponential of $A$ can be expressed as $p(A)$ for a
polynomial $p(z)$.  It is enough to choose $p(z)$ so that it agrees
with the exponential function at the eigenvalues of $A$, and to
sufficiently high order.  Notice in particular that the eigenvalues
of $\exp (A)$ are therefore all exponentials of eigenvalues of $A$.

	We know that the exponential of a linear transformation is
automatically invertible.  Conversely, if $B$ is an invertible linear
transformation on ${\bf C}^n$, is there a linear transformation $A$
on ${\bf C}^n$ such that $\exp (A) = B$?  The answer is yes, and indeed
one can take $A = p(B)$, where $p(z)$ is a polynomial on ${\bf C}$
which satisfies $\exp (p(z)) = z$ at the eigenvalues of $B$, and to
sufficiently high order.

	Now let us consider the real case.  We have seen that the
determinant of the exponential of a linear transformation on ${\bf
R}^n$ is equal to the exponential of the trace of that linear
transformation, and hence is a positive real number.  This is a simple
necessary condition for an invertible linear transformation on ${\bf
R}^n$ to be the exponential of another linear transformation.

	Let $A$ be any linear transformation on ${\bf R}^n$, and let
$\widehat{A}$ denote the unique linear transformation on ${\bf C}^n$
which agree with $A$ on ${\bf R}^n$.  To be more precise,
$\widehat{A}$ is complex-linear, so that $\widehat{A}(i \, v) = i \,
\widehat{A}(v)$, and this ensures that $\widehat{A}$ is determined by
its action on vectors with real coordinates.  Also, $A$ and
$\widehat{A}$ are associated to the same $n \times n$ matrix with real
entries, with respect to the standard bases of ${\bf R}^n$ and ${\bf
C}^n$, respectively.  

	For any polynomial $p(x)$ with real coefficients, we can
define $p(A)$ in the usual manner, and it has the same basic
properties as before.  We can also think of $p$ as a complex
polynomial and consider $p(\widehat{A})$, and it is easy to see that
this is the same as the complex-linear transformation on ${\bf C}^n$
induced by $p(A)$.  In other words,
\begin{equation}
	\widehat{p(A)} = p(\widehat{A}).
\end{equation}

	If $\lambda$ is a real number which is an eigenvalue of $A$,
then $\lambda$ is also an eigenvalue of $\widehat{A}$, using the same
eigenvector in fact, and conversely if $\lambda$ is an eigenvalue of
$\widehat{A}$ which is a real number too, then one can check that
$\lambda$ is an eigenvalue for $A$.  However, in general there can be
complex eigenvalues for $\widehat{A}$, and one can check that if
$\lambda$ is an eigenvalue of $\widehat{A}$, then so is the complex
conjugate $\overline{\lambda}$, and with the same multiplicity as a
zero of the characteristic polynomial $q_{\widehat{A}}$.  Notice
that for $x$ real the characteristic polynomial $q_{\widehat{A}}(x)$
of $\widehat{A}$ is the same as the real version for $A$,
\begin{equation}
	q_{\widehat{A}}(x) = \det (x \, I - A),
\end{equation}
and in particular the two polynomials have the same coefficients,
which are real numbers.

	Suppose that $B$ is an invertible linear transformation on
${\bf R}^n$.  If, for instance, $B$ has no eigenvalues which are
negative real numbers, then there are polynomials $p(z)$ with real
coefficients such that $\exp(p(z)) = z$ to whatever order one might
like at the eigenvalues of $\widehat{B}$.  Consequently, $B = \exp(A)$
with $A = p(B)$, and where $A$ is a linear transformation on ${\bf
R}^n$.

	This is certainly not the whole story however.  Let us mention
two basic examples of linear transformations with positive determinant
and negative real eigenvalues which can and which cannot be
represented as an exponential.  We can look at this in terms of
another correspondence between real and complex-linear
transformations.

	Namely, we can identify ${\bf C}^n$ with ${\bf R}^{2n}$ in the
obvious way, with the real and imaginary parts of the $n$ complex
components of a vector in ${\bf C}^n$ being the $2n$ real components
of the corresponding vector in ${\bf R}^{2n}$.  If $A$ is a linear
transformation on ${\bf C}^n$, let us write $A^\circ$ for the
corresponding real-linear transformation on ${\bf R}^{2n}$.
Notice that
\begin{equation}
	(\alpha_1 \, A_1 + \alpha_2 \, A_2)^\circ
		= \alpha_1 \, A_1^\circ + \alpha_2 \, A_2^\circ
\end{equation}
when $\alpha_1$, $\alpha_2$ are real numbers and $A_1$, $A_2$
are complex-linear transformations on ${\bf C}^n$, and that
$(A_1 \, A_2)^\circ = A_1^\circ \, A_2^\circ$.

	On ${\bf R}^2$, consider the linear transformation $-I$.  This
is a diagonalizable linear transformation with eigenvalue $-1$ of
multiplicity $2$, and the determinant is equal to $1$.  This linear
transformation is the exponential of another linear transformation on
${\bf R}^2$, because one can think of it as a complex-linear
transformation on ${\bf C}$, and convert the realization as an
exponential there to one on ${\bf R}^2$.

	As a different example, suppose that $B$ is a linear
transformation on ${\bf R}^2$ such that the two standard basis vectors
$e_1 = (1, 0)$ and $e_2(0,1)$ are eigenvectors with eigenvalues
$\lambda_1$, $\lambda_2$, respectively, and where $\lambda_1$,
$\lambda_2$ are distinct negative real numbers.  If $B = \exp(A)$ for
some real linear transformation $A$ on ${\bf R}^2$, then $A$, $B$
commute in particular, and it follows that $A(e_1)$, $A(e_2)$ are
eigenvectors for $B$ with eigenvalues $\lambda_1$, $\lambda_2$,
respectively.  Hence $A(e_1)$, $A(e_2)$ should be real multiples
of $e_1$, $e_2$, this leads to a contradiction.

	If $A$ is any complex-linear transformation on ${\bf C}^n$
and $A^\circ$ is the corresponding real-linear transformation on
${\bf R}^{2n}$, then
\begin{equation}
	\det A^\circ = |\det A|^2,
\end{equation}
i.e., the determinant of $A^\circ$ as a real-linear transformation is
equal to the absolute value squared of the determinant of $A$ as a
complex-linear transformation.  This is not too difficult to show,
starting with $n = 1$, for instance.  In particular, $A^\circ$
always has nonnegative determinant.

	There is another well-known simple trick for expressing a
positive power of a linear transformation as a linear combination of
lower powers.  Namely, if $T$ is a linear transformation on ${\bf
R}^n$ or ${\bf C}^n$, then there is a positive integer $k \le n^2$
such that $T^k$ is a linear combination of the identity operator and
$T^j$, $1 \le j < k$, simply because the vector space of linear
transformations on ${\bf R}^n$ or ${\bf C}^n$ has dimension $n^2$.  Of
course the version of this from the Cayley--Hamilton theorem is more
precise and explicit.

	Suppose that $A$ is a linear transformation on ${\bf R}^n$,
and let $\widehat{A}$ denote the corresponding complex-linear
transformation on ${\bf C}^n$.  Of course
\begin{equation}
	\|\widehat{A}\|_{HS} = \|A\|_{HS},
\end{equation}
since $\widehat{A}$ and $A$ correspond to the same $n \times n$
matrix of real numbers, and the norms in question are simply the
square root of the sum of squares of these matrix entries.
Moreover, one can check that
\begin{equation}
	\|\widehat{A}\|_{op} = \|A\|_{op},
\end{equation}
where the left side refers to the operator norm of $\widehat{A}$
as a linear transformation on ${\bf C}^n$, and the right side
refers to the operator norm of $A$ as a linear transformation on
${\bf R}^n$.  

	Also,
\begin{equation}
	\widehat{A}^* = \widehat{(A^*)},
\end{equation}
where the left side is the adjoint of $\widehat{A}$ as a linear
transformation on ${\bf C}^n$, and the right side is the
complex-linear transformation on ${\bf C}^n$ induced by the adjoint of
$A$ as a real-linear transformation on ${\bf R}^n$.  It follows that
if $A$ is an orthogonal linear transformation on ${\bf R}^n$, then
$\widehat{A}$ is a unitary linear transformation.  One can see this as
well using the fact that a real or complex-linear transformation is
orthogonal or unitary, respectively, if and only if it is invertible,
has operator norm equal to $1$, and its inverse has operator norm
equal to $1$.

	If $A$ is an orthogonal or unitary linear transformation
on ${\bf R}^n$ or ${\bf C}^n$, respectively, then
\begin{equation}
	|\det A| = 1.
\end{equation}
Indeed, in this case we have that
\begin{equation}
	1 = \det I = \det (A \, A^*) = (\det A) (\det A^*) = |\det A|^2.
\end{equation}
Alternatively, one can show this using the fact that
\begin{equation}
	|\det A| \le \|A\|_{op}^n.
\end{equation}

	More precisely, if $A$ is a linear transformation on ${\bf C}^n$,
then $\det A$ is the product of the eigenvalues of $A$, according to
their multiplicities as zeros of the characteristic polynomial of $A$,
and it is easy to see that
\begin{equation}
	|\lambda| \le \|A\|_{op}
\end{equation}
for each eigenvalue of $A$.  In the real case one can apply this
argument to the induced complex-linear transformation on ${\bf C}^n$,
which has the same determinant.  As another argument, it is enough
to check that
\begin{equation}
	|\det A| \le 1 \quad\hbox{when}\quad \|A\|_{op} \le 1,
\end{equation}
and for this one can observe that the sequence of linear
transformations $A^k$, $k \ge 1$, is bounded when $\|A\|_{op} \le 1$, and
hence that their determinants are bounded, and hence that the scalars
$(\det A)^k$ are bounded, which implies that $|\det A| \le 1$.

	Suppose that $A$ is a linear transformation on ${\bf R}^n$
or on ${\bf C}^n$ which is anti-self-adjoint, so that
\begin{equation}
	A^* = - A.
\end{equation}
In this case $\exp (A)$ is an orthogonal or unitary transformation, as
appropriate.  The adjoint of $\exp(A)$ is equal to $\exp (A^*)$, which
is the same as $\exp (-A)$ in this case, which is the inverse of $\exp
(A)$.

	Let us consider next the question of when an orthogonal linear
transformation on ${\bf R}^n$ or a unitary transformation on ${\bf
C}^n$ can be expressed as the exponential of a self-adjoint linear
transformation.  To do this we digress a bit for some general matters
about linear transformations.  We begin with the complex case.

	A linear transformation $T$ on ${\bf C}^n$ is said to be
\emph{normal} if $T$ commutes with its adjoint, which is to say that
\begin{equation}
	T^* \, T = T \, T^*.
\end{equation}
We can write any linear transformation $T$ on ${\bf C}^n$ as $T_1 + i
\, T_2$, where $T_1$, $T_2$ are the self-adjoint linear
transformations given by
\begin{equation}
	T_1 = \frac{1}{2} (T + T^*), \quad
		T_2 = \frac{1}{2i} (T - T^*),
\end{equation}
and the condition of normality is equivalent to saying that $T_1$,
$T_2$ commute.  Note that unitary transformations are normal.

	We already know that if $B$ is a self-adjoint linear
transformation, then there is an orthonormal basis of the underlying
vector space consisting of eigenvectors of $B$.  Given two
self-adjoint linear transformations which commute, one can find an
orthonormal basis consisting of vectors which are eigenvectors for
both linear transformations.  Conversely, for a fixed basis, any two
linear transformations for which vectors in the basis are eigenvectors
clearly commute with each other.

	As a result, if $T$ is a normal linear transformation on ${\bf
C}^n$, then there is an orthonormal basis of ${\bf C}^n$ consisting of
eigenvectors of $T$.  In particular this applies to unitary
transformations, for which the corresponding eigenvalues are complex
numbers with modulus $1$.  As a result, if $U$ is a unitary linear
transformation on ${\bf C}^n$, then there is an anti-self-adjoint
linear transformation $A$ on ${\bf C}^n$ such that $\exp (A) = U$,
and indeed one can take $A$ to be diagonalized by the same basis
as for $U$, with imaginary eigenvalues.

	Now let us consider the real case.  For this we cannot use
the trick of writing an anti-self-adjoint linear transformation
as ``$i$'' times a self-adjoint linear transformation.  There are
other things that we can do, however.

	Thus we let $A$ be an anti-self-adjoint linear transformation
on ${\bf R}^n$.  Notice that
\begin{equation}
	\langle A(v), v \rangle = 0
\end{equation}
for all vectors $v \in {\bf R}^n$, and that in particular a nonzero
vector $v$ in ${\bf R}^n$ is an eigenvector for $A$ only if the
corresponding eigenvalue is equal to $0$, so that $v$ lies in the
kernel of $A$.  Also, for each vector $w$ in ${\bf R}^n$, $A(w)$ is
orthogonal to every vector in the kernel of $A$.

	A basic trick to study an anti-self-adjoint linear
transformation $A$ is to consider $A^2$, which is self-adjoint and has
the same kernel as $A$ does.  If $v$ is a vector in ${\bf R}^n$ which
is an eigenvector $A^2$ with eigenvalue $\lambda$, then $A(v)$ is an
eigenvector for $A^2$ with eigenvalue $\lambda$ too, and of course
$A^2(v)$ is a multiple of $v$.  As a result, if $\lambda$ is a
negative real number which is an eigenvalue of $A^2$, then one can
check that the corresponding eigenspace
\begin{equation}
	\{ v \in {\bf R}^n : A^2(v) = \lambda \, v \}
\end{equation}
has \emph{even} dimension.

	If $T$ is an orthogonal linear transformation on ${\bf R}^n$,
then we can write $T = T_1 + T_2$, where $T_1 = (T + T^*)/2$
is self-adjoint, $T_2 = (T - T^*)/2$ is anti-self-adjoint, 
$T_1$, $T_2$ commute, and
\begin{equation}
	\langle T_1(v), T_2(v) \rangle = 0
\end{equation}
for all $v \in {\bf R}^n$.  If $\lambda_1$, $\lambda_2$ are
eigenvalues of $T_1$, $T_2^2$ such that the joint eigenspace
\begin{equation}
	\{ v \in {\bf R}^n : T_1(v) = \lambda_1 \, v,
				\ T_2^2(v) = \lambda_2 \, v \}
\end{equation}
is nontrivial, then either $\lambda_2 = 0$ and $\lambda_1 = \pm 1$, or
$\lambda_2 < 0$, $\lambda_1^2 - \lambda_2 = 1$, and the joint
eigenspace has even dimension.  One can show that the parity of the
number of times that $\lambda_1 = -1$ is even or odd according to
whether $\det T$ is $1$ or $-1$, and that when $\det T = 1$, 
there is an anti-self-adjoint linear transformation $A$ on ${\bf R}^n$
such that $T = \exp(A)$.

	If $A$ is a self-adjoint linear transformation on ${\bf R}^n$
or on ${\bf C}^n$, then $\exp (A)$ is also self-adjoint, and in fact
$\exp (A)$ is positive-definite, because $\exp (A) = B^2$ where $B$ is
the self-adjoint linear transformation $\exp (A/2)$.  Conversely,
every self-adjoint positive-definition linear transformation $P$ on
${\bf R}^n$ or on ${\bf C}^n$ can be realized as $\exp (A)$ for a
self-adjoint linear transformation $A$.  This can be seen using an
orthonormal basis of eigenvectors for $P$, and indeed one can choose
$A$ so that the vectors in the same basis are eigenvectors for $A$.

	Actually, if $A$ is a self-adjoint linear transformation on
${\bf R}^n$ or on ${\bf C}^n$, then there is an orthonormal basis of
eigenvectors for $A$, and of course these same vectors are
eigenvectors for $\exp (A)$, with the eigenvalues for $\exp (A)$ being
the exponentials of the corresponding eigenvalues for $A$.  As a
result, for a given self-adjoint positive-definite linear
transformation $P$ on ${\bf R}^n$ or ${\bf C}^n$, the self-adjoint
linear transformation $A$ on ${\bf R}^n$ or ${\bf C}^n$, respectively,
such that $\exp (A) = P$ is \emph{unique}.  This is analogous to the
situation for the ordinary exponential function on real numbers,
while in the complex case one can have different numbers or linear
transformations whose exponentials are equal to each other.

	There is a natural mapping from the group of invertible
linear transformations on ${\bf R}^n$ or on ${\bf C}^n$ onto the
self-adjoint positive-definite linear transformations on the same
space, given by
\begin{equation}
	T \mapsto T \, T^*.
\end{equation}
If $T$ is an invertible linear transformation on ${\bf R}^n$ or on
${\bf C}^n$ and $R$ is an orthogonal or unitary linear transformation
on the same space, as appropriate, then $T$ and $T \, R$ are sent by
the mapping just defined to the same positive-definite linear
transformation, since
\begin{equation}
	(T \, R) (T \, R)^* = T \, R \, R^* \, T^* = T \, T^*.
\end{equation}
Conversely, if $T$, $T'$ are invertible linear transformations on
${\bf R}^n$ or on ${\bf C}^n$ such that $T' \, (T')^* = T \, T^*$,
then there is an orthogonal or unitary linear transformation $R$,
as appropriate, such that
\begin{equation}
	T' = T \, R.
\end{equation}

	Also, every self-adjoint positive-definite linear
transformation $P$ on ${\bf R}^n$ or on ${\bf C}^n$ arises this
manner, and in fact can be written in a unique manner as $Q^2$ for a
self-adjoint positive-definite linear transformation $Q$.  If $A$ is
an invertible linear mapping on ${\bf R}^n$ or on ${\bf C}^n$, then we
get a nice action of $A$ on the self-adjoint positive-definite linear
transformations by the formula
\begin{equation}
	P \mapsto A \, P \, A^*.
\end{equation}
If $A_1$ and $A_2$ are two invertible linear transformations on ${\bf
R}^n$ or on ${\bf C}^n$ such that $A_1 \, P \, A_1^* = A_2 \, P \,
A_2^*$ for all self-adjoint positive-definite linear transformations
$P$, then $A_1 = A_2$, and if $P_1$, $P_2$ are two self-adjoint
positive-definite linear transformations on ${\bf R}^n$ or on ${\bf C}^n$,
then there is an invertible linear transformation $A$ on the same space
such that $P_2 = A \, P_1 \, A^*$.

\section{Spaces of matrices}
\label{spaces of matrices}
\setcounter{equation}{0}

	As before, we write $\mathcal{L}({\bf R}^n)$,
$\mathcal{L}({\bf C}^n)$ for the spaces of real and complex-linear
mappings on ${\bf R}^n$ and ${\bf C}^n$, respectively.  We also write
$GL({\bf R}^n)$, $GL({\bf C}^n)$ for the \emph{general linear groups}
of invertible real and complex-linear transformations on ${\bf R}^n$,
${\bf C}^n$, respectively.  We can identify $\mathcal{L}({\bf R}^n)$,
$\mathcal{L}({\bf C}^n)$ with ${\bf R}^{n^2}$, ${\bf C}^{n^2}$ using
the standard correspondence between linear transformations and
matrices, and in this way we have that $GL({\bf R}^n)$, $GL({\bf C}^n)$
are open subsets of $\mathcal{L}({\bf R}^n)$, $\mathcal{L}({\bf C}^n)$,
respectively.

	The determinant can be viewed as a homogeneous polynomial of
degree $n$ on $\mathcal{L}({\bf R}^n)$, $\mathcal{L}({\bf C}^n)$, and
$GL({\bf R}^n)$, $GL({\bf C}^n)$ can be described as the subsets of
$\mathcal{L}({\bf R}^n)$, $\mathcal{L}({\bf C}^n)$ defined by the
condition
\begin{equation}
	\det T \ne 0.
\end{equation}
At the identity operator, the differential of the determinant can
be identified with the trace, since
\begin{equation}
	\biggl(\frac{d}{dr} \det(I + r \, A) \biggr)_{r = 0}
		= \tr A
\end{equation}
for any linear transformation $A$ on ${\bf R}^n$ or ${\bf C}^n$.
If $T$ is an invertible linear transformation on ${\bf R}^n$
or on ${\bf C}^n$ and $A$ is another linear transformation on
the same space, then the differential of the determinant at $T$
in the direction $A$ can be expressed as
\begin{equation}
	d(\det)_T(A) = (\det T) \tr (T^{-1} \, A),
\end{equation}
since
\begin{eqnarray}
  \biggl(\frac{d}{dr} \det (T + r \, A)\biggr)_{r = 0}
	& = & (\det T) 
	   \biggl(\frac{d}{dr} \det (I + r \, T^{-1} \, A)\biggr)_{r = 0}
								\\
	& = & (\det T) \tr (T^{-1} \, A).
\end{eqnarray}

	We write $SL({\bf R}^n)$, $SL({\bf C}^n)$ for the subgroups of
$GL({\bf R}^n)$, $GL({\bf C}^n)$ consisting of linear transformations
with determinant equal to $1$.  These are nice submanifolds of
$GL({\bf R}^n)$, $GL({\bf C}^n)$, because the differential of the 
determinant is not equal to $0$ at any point in $SL({\bf R}^n)$,
$SL({\bf C}^n)$, or in $GL({\bf R}^n)$, $GL({\bf C}^n)$, for that
matter.  Also we have the maps
\begin{equation}
	T \mapsto (\det T)^{-1/n} \, T
\end{equation}
from invertible linear transformations on ${\bf R}^n$, ${\bf C}^n$
to linear transformations with determinant $1$, at least if we
restrict our attention to $T$'s with $\det T > 0$ in the real case
and $T$'s with $\det T$ in a nice region in ${\bf C}$ that contains
$1$ and on which $z^{-1/n}$ can be defined in the complex case.

	If $T$ is an element of $GL({\bf R}^n)$ or $GL({\bf C}^n)$,
then the space of tangent vectors to $GL({\bf R}^n)$ or $GL({\bf
C}^n)$ at $T$, as appropriate, can be identified with
$\mathcal{L}({\bf R}^n)$ or $\mathcal{L}({\bf C}^n)$, as appropriate,
since $GL({\bf R}^n)$, $GL({\bf C}^n)$ are open subsets of
$\mathcal{L}({\bf R}^n)$, $\mathcal{L}({\bf C}^n)$, respectively.  If
$T$ is an element of $SL({\bf R}^n)$ or $SL({\bf C}^n)$, then the
space of tangent vectors to $SL({\bf R}^n)$ or $SL({\bf C}^n)$
at $T$ is equal to the space of $A$ in $\mathcal{L}({\bf R}^n)$
or $\mathcal{L}({\bf C}^n)$ such that
\begin{equation}
	\tr (T^{-1} \, A) = 0,
\end{equation}
as appropriate.  In other words, these are the tangent vectors
to $GL({\bf R}^n)$ or $GL({\bf C}^n)$ at $T$, respectively,
which lie in the kernel of the differential of the determinant
function at $T$.

	An important mapping on $GL({\bf R}^n)$ or $GL({\bf C}^n)$ is
the one defined by
\begin{equation}
	F(T) = T^{-1},
\end{equation}
and which also sends $SL({\bf R}^n)$ or $SL({\bf C}^n)$ to itself, as
appropriate.  For each invertible linear transformation $T$ the
differential of this mapping at $T$ can be expressed as
\begin{equation}
	dF_T(A) = - T^{-1} \, A \, T^{-1},
\end{equation}
which is to say that
\begin{equation}
	\biggl(\frac{d}{dr} (T + r \, A)^{-1} \biggr)_{r = 0}
		= - T^{-1} \, A \, T^{-1}.
\end{equation}
Indeed,
\begin{eqnarray}
	(T + r \, A)^{-1} 
		& = & (I + r \, T^{-1} \, A)^{-1} \, T^{-1}	\\
		& = & (I - r \, T^{-1} \, A + O(r^2)) \, T^{-1}
			= T^{-1} - r \, T^{-1} \, A \, T^{-1} + O(r^2).
							\nonumber
\end{eqnarray}

	If $T$ is an invertible linear transformation on ${\bf R}^n$
and $A$, $B$ are linear transformations on ${\bf R}^n$, then set
\begin{equation}
	\langle A, B \rangle_T = \tr (T^{-1} \, A \, T^{-1} \, B).
\end{equation}
This is a symmetric bilinear form in $A$, $B$ for each $T$,
which is to say that $\langle A, B \rangle_T$ is a linear function
of $A$ for each $B$ and $T$, a linear function of $B$ for each $A$
and $T$, and that in fact
\begin{equation}
	\langle A, B \rangle_T = \langle B, A \rangle_T
\end{equation}
for all $A$, $B$, $T$, so that linearity in $A$ and $B$ are
equivalent to each other.  Moreover, this bilinear form is
nondegerate, which means that for each $T$ and for each $A \ne 0$
there is a $B$ such that $\langle A, B \rangle_T \ne 0$.

	In the complex case, if $T$ is an invertible linear
transformation on ${\bf C}^n$ and $A$, $B$ are linear transformations
on ${\bf C}^n$, then
\begin{equation}
	\tr (T^{-1} \, A \, T^{-1} \, B)
\end{equation}
is a complex-valued symmetric bilinear form in $A$, $B$ for each $T$
which is nondegenerate.  To get a real-valued quantity, we set
\begin{equation}
	\langle A, B \rangle_T = \re \tr(T^{-1} \, A \, T^{-1} \, B),
\end{equation}
i.e., we take the real part of the trace.  This is still
real-bilinear, which means that it is real-linear in each of $A$ and
$B$, and symmetric and nondegenerate.

	Notice that $\langle A, B \rangle_T$ depends smoothly on $T$
for $T$ in $GL({\bf R}^n)$ or $GL({\bf C}^n)$, as appropriate.  As a
result, $\langle A, B \rangle_T$ is said to define a semi-Riemannian
structure, also known as a pseudo-Riemannian structure or a Riemannian
structure with signature, on $GL({\bf R}^n)$ or $GL({\bf C}^n)$, as
appropriate.  On $GL({\bf C}^n)$, if we did not take the real part,
then we would have a holomorphic semi-Riemannian structure.

	At points $T$ in $SL({\bf R}^n)$ or $SL({\bf C}^n)$, we can
restrict our attention to $A$, $B$ which are in the tangent space of
$SL({\bf R}^n)$ or $SL({\bf C}^n)$ at $T$, as appropriate.
Explicitly, this means that we restrict our attention to $A$, $B$ such
that
\begin{equation}
	\tr (T^{-1} \, A) = \tr (T^{-1} \, B) = 0.
\end{equation}
This leads to semi-Riemannian structures on $SL({\bf R}^n)$, $SL({\bf
C}^n)$, and a holomorphic semi-Riemannian structure on $SL({\bf C}^n)$
if we do not take the real part.

	For each linear transformation $Z$ on ${\bf R}^n$ or on ${\bf
C}^n$, define linear transformations $\lambda_Z$, $\rho_Z$ on
$\mathcal{L}({\bf R}^n)$ or on $\mathcal{L}({\bf C}^n)$, respectively,
by
\begin{equation}
	\lambda_Z(T) = Z \, T, \quad \rho_Z(T) = T \, Z,
\end{equation}
which is to say that $\lambda_Z$, $\rho_Z$ correspond to left and
right multiplication by $Z$.  These are linear transformations, and
in particular their differentials are given by themselves,
\begin{equation}
	(d \lambda_Z)_T(A) = \lambda_Z(A),
		(d \rho_Z)_T(A) = \rho_Z(A)
\end{equation}
for all linear transformations $T$, $A$ on ${\bf R}^n$ or on ${\bf
C}^n$, as appropriate.  If $Z$ is an invertible linear transformation
on ${\bf R}^n$ or on ${\bf C}^n$, then $\lambda_Z$, $\rho_Z$ map
$GL({\bf R}^n)$ or $GL({\bf C}^n)$ onto itself, as appropriate, and if
\begin{equation}
	\det Z = 1,
\end{equation}
then $\lambda_Z$, $\rho_Z$ map $SL({\bf R}^n)$ or $SL({\bf C}^n)$
onto itself, as appropriate.

	If $Z$ is an invertible linear transformation on ${\bf R}^n$
or on ${\bf C}^n$, then the mappings $\lambda_Z$, $\rho_Z$ on $GL({\bf
R}^n)$ or on $GL({\bf C}^n)$ preserve the semi-Riemannian structure
$\langle \cdot, \cdot \rangle_T$.  In other words, if $T$ is an
element of $GL({\bf R}^n)$ or $GL({\bf C}^n)$ and $A$, $B$ are
elements of $\mathcal{L}({\bf R}^n)$ or $\mathcal{L}({\bf C}^n)$, as
appropriate, which we view as tangent vectors to $GL({\bf R}^n)$
or $GL({\bf C}^n)$ at $T$, then
\begin{equation}
	\langle (d \lambda_Z)_T(A), (d \lambda_Z)_T(B) \rangle_{\lambda_Z(T)}
		= \langle A, B \rangle_T
\end{equation}
and
\begin{equation}
	\langle (d \rho_Z)_T(A), (d \rho_Z)_T(B) \rangle_{\rho_Z(T)}
		= \langle A, B \rangle_T.
\end{equation}
This is easy to verify.  By restriction, if $\det Z = 1$,
so that $\lambda_Z$, $\rho_Z$ can be viewed as defining mappings
on $SL({\bf R}^n)$ or $SL({\bf C}^n)$, as appropriate, then
$\lambda_Z$, $\rho_Z$ preserve the restriction of our semi-Riemannian
structures to $SL({\bf R}^n)$ or $SL({\bf C}^n)$.

	One can also check that $F(T) = T^{-1}$ preserves the
semi-Riemannian structures on $GL({\bf R}^n)$, $GL({\bf C}^n)$.  That
is, if $T$ is an element of $GL({\bf R}^n)$ or $GL({\bf C}^n)$, and if
$A$, $B$ are elements of $\mathcal{L}({\bf R}^n)$ or $\mathcal{L}({\bf
C}^n)$, respectively, which we can view as tangent vectors to $GL({\bf
R}^n)$ or $GL({\bf C}^n)$ at $T$, then
\begin{equation}
	\langle (d F)_T(A), (d F)_T(B) \rangle_{F(T)}
		= \langle A, B \rangle_T.
\end{equation}
Also, $SL({\bf R}^n)$ and $SL({\bf C}^n)$ are invariant under $F(T) =
T^{-1}$, and the restriction of the semi-Riemannian structure to
$SL({\bf R}^n)$, $SL({\bf C}^n)$ is preserved by $F$, since this holds
on $GL({\bf R}^n)$, $GL({\bf C}^n)$.

	In the complex case, let us note that $\lambda_Z$, $\rho_Z$
are complex-linear transformations, and $F(T) = T^{-1}$ is a
holomorphic transformation, and that they preserve the holomorphic
version of the semi-Riemannian structure on $GL({\bf C}^n)$
and $SL({\bf C}^n)$.

	Of course we can define a flat semi-Riemannian metric on
$\mathcal{L}({\bf R}^n)$ by saying that if $T$, $A$, $B$ are linear
transformations on ${\bf R}^n$, where we think of $A$, $B$ as tangent
vectors to $\mathcal{L}({\bf R}^n)$ at $T$, then the inner product of
these two tangent vectors associated to $T$ is given by
\begin{equation}
\label{tr (A B)}
	\tr (A \, B).
\end{equation}
In the complex case the same formula defines a holomorphic
semi-Riemannian structure on $\mathcal{L}({\bf C}^n)$, and to get an
ordinary semi-Riemannian structure one should take the real part.
That (\ref{tr (A B)}) does not depend on $T$ reflects the fact that
these semi-Riemannian structures are flat.  Of course we can restrict
these semi-Riemannian structures to the subspaces $\mathcal{L}_0({\bf
R}^n)$, $\mathcal{L}_0({\bf C}^n)$ of linear transformations with
trace equal to $0$.

	The exponential mapping defines a mapping from
$\mathcal{L}({\bf R}^n)$, $\mathcal{L}({\bf C}^n)$ to $GL({\bf R}^n)$,
$GL({\bf C}^n)$ respectively, sending $0$ to $I$ and with
\begin{equation}
	d \exp_0 (A) = A.
\end{equation}
In particular, the standard flat metric at $0$ on $\mathcal{L}({\bf R}^n)$,
$\mathcal{L}({\bf C}^n)$ corresponds exactly to the semi-Riemannian
structure on $GL({\bf R}^n)$, $GL({\bf C}^n)$ at $I$, respectively,
under the differential of the exponential mapping at $0$.  In fact,
\begin{eqnarray}
\label{langle d exp_T(A), d exp_T(B) rangle_{exp T}}
\lefteqn{\langle d \exp_T(A), d \exp_T(B) \rangle_{\exp T}}	\\
	   & & 	= \tr (\exp (-T)) (d \exp_T(A)) (\exp (-T)) (d \exp_T (B))
							\nonumber
\end{eqnarray}
agrees with
\begin{equation}
	\tr A \, B
\end{equation}
to another term in the Taylor expansion, which is to say up to terms
of order $O(\|T\|_{op}^2)$.  In other words, using $\exp (-T) = I - T
+O(\|T\|_{op}^2)$,
\begin{equation}
	d \exp_T (A) = A + \frac{1}{2}(T \, A + A \, T) + O(\|T\|_{op}^2),
\end{equation}
and similarly for $B$, one can check that the terms in (\ref{langle d
exp_T(A), d exp_T(B) rangle_{exp T}}) with no $T$'s reduce to $\tr (A
\, B)$, and that the terms with exactly one $T$ cancel out.  In the
complex case, let us note that the exponential mapping is holomorphic,
and one has the analogous statement about the holomorphic
semi-Riemannian metrics on $\mathcal{L}({\bf C}^n)$ and $GL({\bf
C}^n)$ agreeing at $0$ up to terms of order $O(\|T\|_{op}^2)$.

	As a consequence, for each linear transformation $A$ on ${\bf
R}^n$ or on ${\bf C}^n$, $\exp (t \, A)$ satisfies the equation for
geodesics at $t = 0$.  This is because $t \, A$ is simply a straight
line in $\mathcal{L}({\bf R}^n)$ or $\mathcal{L}({\bf C}^n)$, as
appropriate, and thus satisfies the equation for geodesics there with
respect to the flat semi-Riemannian structure being used, and because
the exponential mapping takes the flat semi-Riemannian structures on
$\mathcal{L}({\bf R}^n)$, $\mathcal{L}({\bf C}^n)$ around $0$ to the
semi-Riemannian structures on $GL({\bf R}^n)$, $GL({\bf C}^n)$ around
$I$ to sufficient precision, as in the preceding paragraph.  In fact,
it follows that $\exp (t \, A)$ satisfies the equation for geodesics
for all $t$, because one can use the invariance of the semi-Riemannian
metrics on $\mathcal{L}({\bf R}^n)$, $\mathcal{L}({\bf C}^n)$ under
ordinary translations and the invariance of the semi-Riemannian
structures on $GL({\bf R}^n)$, $GL({\bf C}^n)$ under left and right
multiplication by invertible linear transformations to reduce the
case of general $t$ to $t = 0$.  In the complex case, one can take
$t$ to be a complex parameter, and say that $\exp (t \, A)$ is a 
holomorphic geodesic in $GL({\bf C}^n)$ with respect to the
holomorphic semi-Riemannian structure as before.

	If $A$, $Y$ are linear transformations on ${\bf R}^n$ or on
${\bf C}^n$ with $Y$ invertible, then $Y \, exp (t \, A)$ defines a
geodesic in $GL({\bf R}^n)$ or $GL({\bf C}^n)$, as appropriate, again
using the invariance of the semi-Riemannian structures that we have
defined.  This is equivalent to saying that if $B$, $Y$ are linear
transformations on ${\bf R}^n$ or on ${\bf C}^n$, then $\exp (t \, B)
\, Y$ defines a geodesic in $GL({\bf R}^n)$ or $GL({\bf C}^n)$, as
appropriate, since $Y \, \exp (t \, A) = \exp (t \, B) \, Y$ with $B =
Y \, A \, Y^{-1}$, and anyway our semi-Riemannian structures on the
general linear groups are invariant under both left and right
multiplications.  This accounts for all of the geodesics, because the
equation for geodesics are described by a second-order differential
equation, and thus a geodesic is characterized by a point that it
passes through and the tangent vector corresponding to its derivative
at that point.

	The preceding discussion can also be applied to the
restriction of the exponential mapping to the subspaces
$\mathcal{L}_0({\bf R}^n)$, $\mathcal{L}_0({\bf C}^n)$ of
$\mathcal{L}({\bf R}^n)$, $\mathcal{L}({\bf C}^n)$ taking values in
the subgroups $SL({\bf R}^n)$, $SL({\bf C}^n)$ of $GL({\bf R}^n)$,
$GL({\bf C}^n)$, using the restrictions of the corresponding
semi-Riemannian metrics.  In particular, $SL({\bf R}^n)$, $SL({\bf
C}^n)$ are \emph{totally geodesic} submanifolds of $GL({\bf R}^n)$,
$GL({\bf C}^n)$.  That is, a geodesic in $GL({\bf R}^n)$ or $GL({\bf
C}^n)$ which passes through $SL({\bf R}^n)$ or $SL({\bf C}^n)$ and is
tangent to the special linear group at the point of intersection stays
in the special linear group.  Note that $SL({\bf C}^n)$ is a complex
submanifold of $GL({\bf C}^n)$.

	Fix a positive integer $n$, and let $\mathcal{F}$ be a
\emph{flag} in ${\bf R}^n$ or ${\bf C}^n$, which is to say a
family $L_1, L_2, \ldots, L_k$ of distinct nontrivial proper
linear subspaces of ${\bf R}^n$ or ${\bf C}^n$ with
\begin{equation}
	L_1 \subseteq L_2 \subseteq \cdots \subseteq L_k.
\end{equation}
Thus $k$ is a positive integer strictly less than $n$, called
the \emph{length} of the flag.  It may be that $k = 1$, so that
the flag consists of a single nontrivial proper linear subspace.

	If $\mathcal{F}$ is a flag in ${\bf R}^n$ or in ${\bf C}^n$,
then we write $\mathcal{L}_\mathcal{F}({\bf R}^n)$ or
$\mathcal{L}_\mathcal{F}({\bf C}^n)$ for the space of linear
transformations $A$ on ${\bf R}^n$ or ${\bf C}^n$, as appropriate,
such that $A(L_j) \subseteq L_j$ for each of the linear subspaces
$L_j$ in the flag.  Thus $\mathcal{L}_\mathcal{F}({\bf R}^n)$,
$\mathcal{L}_\mathcal{F}({\bf C}^n)$ are themselves linear subspaces
of $\mathcal{L}({\bf R}^n)$, $\mathcal{L}({\bf C}^n)$, respectively,
which are also closed under taking products of linear transformations.
By using bases for ${\bf R}^n$ or ${\bf C}^n$, as appropriate, which
are suitably adapted to the flag $\mathcal{F}$, one can also
characterize the linear transformations in
$\mathcal{L}_\mathcal{F}({\bf R}^n)$ or $\mathcal{L}_\mathcal{F}({\bf
C}^n)$, as appropriate, in terms of matrices with certain entries
equal to $0$.

	Similarly, if $\mathcal{F}$ is a flag in ${\bf R}^n$ or in
${\bf C}^n$, then we write $GL_\mathcal{F}({\bf R}^n)$ or
$GL_\mathcal{F}({\bf C}^n)$ for the space of invertible linear
transformations $T$ on ${\bf R}^n$ or ${\bf C}^n$, as appropriate,
such that $T(L_j) = L_j$ for each linear subspace $L_j$ in the flag.
This is equivalent to saying that
\begin{equation}
	GL_\mathcal{F}({\bf R}^n) 
		= GL({\bf R}^n) \cap \mathcal{L}_\mathcal{F}({\bf R}^n)
\end{equation}
and
\begin{equation}
	GL_\mathcal{F}({\bf C}^n)
		= GL({\bf C}^n) \cap \mathcal{L}_\mathcal{F}({\bf C}^n).
\end{equation}
Furthermore, let us put
\begin{equation}
	SL_\mathcal{F}({\bf R}^n)
		= SL({\bf R}^n) \cap \mathcal{L}_\mathcal{F}({\bf R}^n)
\end{equation}
and
\begin{equation}
	SL_\mathcal{F}({\bf C}^n)
		= SL({\bf C}^n) \cap \mathcal{L}_\mathcal{F}({\bf C}^n).
\end{equation}

	The exponential mapping can be restricted to
$\mathcal{L}_\mathcal{F}({\bf R}^n)$ or $\mathcal{L}_\mathcal{F}({\bf
C}^n)$ to get a mapping into $GL_\mathcal{F}({\bf R}^n)$ or
$GL_\mathcal{F}({\bf C}^n)$, as appropriate.  One can restrict a bit
further to linear transformations in $\mathcal{L}_\mathcal{F}({\bf
R}^n)$ or $\mathcal{L}_\mathcal{F}({\bf C}^n)$ with trace equal to
$0$, which the exponential mapping sends to linear transformations in
$SL_\mathcal{F}({\bf R}^n)$ or $SL_\mathcal{F}({\bf C}^n)$.  As
before, one can account for all of the geodesics in
$GL_\mathcal{F}({\bf R}^n)$, $GL_\mathcal{F}({\bf C}^n)$ or in
$SL_\mathcal{F}({\bf R}^n)$, $SL_\mathcal{F}({\bf C}^n)$ using
exponentials, and these define totally geodesic submanifolds of
$GL({\bf R}^n)$, $GL({\bf C}^n)$, respectively.

	Let us write $\mathcal{S}({bf R}^n)$, $\mathcal{S}({\bf C}^n)$
for the real vector spaces of self-adjoint linear transformations on
${\bf R}^n$, ${\bf C}^n$, respectively, and $\mathcal{S}_+({\bf
R}^n)$, $\mathcal{S}_+({\bf C}^n)$ for their open cones of
positive-definite linear transformations.  These subsets are invariant
under the transformation $F(T) = T^{-1}$, and also under the action
\begin{equation}
	T \mapsto Z^* \, T \, Z
\end{equation}
for each $Z$ in $GL({\bf R}^n)$, $GL({\bf C}^n)$, as appropriate.  In
fact, this action is \emph{transitive}, which is to say that for each
$T_1$, $T_2$ in $\mathcal{S}_+({\bf R}^n)$ or in $\mathcal{S}_+({\bf
C}^n)$ there is a $Z$ in $GL({\bf R}^n)$ or $GL({\bf C}^n)$, as
appropriate, such that $Z^* \, T_1 \, Z = T_2$.

	The restriction of our semi-Riemannian structures on $GL({\bf
R}^n)$, $GL({\bf C}^n)$ to $\mathcal{S}_+({\bf R}^n)$,
$\mathcal{S}_+({\bf C}^n)$, respectively, are \emph{Riemannian}
metrics, which is to say that they are positive definite.  Of course
these Riemannian metrics are invariant under the transformations
preserving $\mathcal{S}_+({\bf R}^n)$, $\mathcal{S}_+({\bf C}^n)$
mentioned in the previous paragraph.  The exponential mapping sends
$\mathcal{S}({\bf R}^n)$, $\mathcal{S}({\bf R}^n)$ onto
$\mathcal{S}_+({\bf R}^n)$, $\mathcal{S}_+({\bf C}^n)$, respectively,
and the geodesics in $\mathcal{S}_+({\bf R}^n)$, $\mathcal{S}_+({\bf
C}^n)$ through $I$ are exactly the curves $\exp (t \, A)$, with $A$ in
$\mathcal{S}({\bf R}^n)$, $\mathcal{S}({\bf C}^n)$, respectively.  The
geodesics through a point $T = Z^* \, Z$ are of the form $Z^* \,
\exp(t \, A) \, Z$.

	The orthogonal and unitary groups on ${\bf R}^n$, ${\bf C}^n$
are denoted $O({\bf R}^n)$ and $U({\bf C}^n)$ and consist of the
invertible linear transformations $T$ which are orthogonal or unitary,
respectively, which is to say that
\begin{equation}
	T^* \, T = T \, T^* = I.
\end{equation}
It is enough to have $T^* \, T = I$, and at a point $T$ in $O({\bf
R}^n)$ or $U({\bf C}^n)$ the tangent space to the orthogonal or
unitary group consists of the linear transformations $A$ on ${\bf
R}^n$ or ${\bf C}^n$, as appropriate, such that
\begin{equation}
	T^* \, A + A^* \, T = 0.
\end{equation}
Let us note that the orthogonal and unitary groups are compact smooth
submanifolds of the vector spaces of all linear transformations on
${\bf R}^n$, ${\bf C}^n$, and of $GL({\bf R}^n)$, $GL({\bf C}^n)$
in particular.

	Again we can restrict our semi-Riemannian structures from
$GL({\bf R}^n)$ or $GL({\bf C}^n)$ to the submanifolds given by the
orthogonal and unitary groups, respectively.  Now these restricted
structures are negative-definite, so that their negatives are
Riemannian metrics.  Using the group structure we again have the
mappings $\lambda_Z(T) = Z \, T$ and $\rho_Z(T) = T \, Z$ which send
the orthogonal and unitary groups to themselves as long as $Z$ also
lies in the orthogonal or unitary group, and also the mapping $F(T) =
T^{-1}$ takes the orthogonal and unitary groups to themselves as well.
The negative Riemannian metrics on the orthogonal and unitary groups
are preserves by these transformations.  If $A$ is an
anti-self-adjoint linear transformation on ${\bf R}^n$ or on ${\bf
C}^n$, then $\exp (t \, A)$ defines a geodesic in $O({\bf R}^n)$ or
$U({\bf C}^n)$, as appropriate, and this accounts for all geodesics in
the orthogonal and unitary groups through $I$, and hence for all
geodesics if one also takes into account the left or right translation
mappings $\lambda_Z$, $\rho_Z$, as before.

	In these various case one can restrict further to linear
transformations with determinant equal to $1$.  Let us write
$\mathcal{M}({\bf R}^n)$, $\mathcal{M}({\bf C}^n)$ for the
hypersurfaces in $\mathcal{S}_+({\bf R}^n)$, $\mathcal{S}_+({\bf
C}^n)$ consisting of linear transformations with determinant equal to
$1$, and $SO({\bf R}^n)$ and $SU({\bf C}^n)$ for the special
orthogonal and unitary groups on ${\bf R}^n$ and ${\bf C}^n$, which
are the subgroups of $O({\bf R}^n)$, $U({\bf C}^n)$ determined by the
condition that the determinant of the corresponding linear
transformation be equal to $1$.  There are similar considerations as
before concerning tangent vectors, Riemannian structures, geodesics,
and so on.

\section{Some geometric situations}
\label{some geometric situations}
\setcounter{equation}{0}

	As usual, ${\bf Z}$ denotes the integers, and ${\bf Z}^n$
consists of $n$-tuples of integers.  Sometimes we might refer to ${\bf
Z}^n$ as the \emph{standard integer lattice} in ${\bf R}^n$.  If we
say that $L$ is a \emph{lattice} in ${\bf R}^n$, then we mean that
there is an invertible linear transformation $A$ on ${\bf R}^n$ such
that
\begin{equation}
	L = A({\bf Z}^n).
\end{equation}

	If $L$ is a lattice in ${\bf R}^n$, then we can form the
quotient space ${\bf R}^n / L$.  That is, two vectors $x$, $y$ in
${\bf R}^n$ are identified in the quotient if their difference $x - y$
lies in $L$.  In particular, we get a canonical quotient mapping
\begin{equation}
	q : {\bf R}^n \to {\bf R}^n / L
\end{equation}
which sends a vector $x$ in ${\bf R}^n$ to the corresponding
element of the quotient.

	Now, with respect to ordinary vector addition, ${\bf R}^n$ is
an abelian group, and a lattice $L$ is a subgroup of ${\bf R}^n$.  We
can think of the quotient ${\bf R}^n / L$ as a quotient in the sense
of group theory.  The quotient is an abelian group under addition, and
the canonical quotient mapping is a group homomorphism.

	We can also look at the quotient ${\bf R}^n / L$ in terms of
topology.  Namely, it inherits a topology from the one on ${\bf R}^n$
so that the canonical quotient mapping is an open continuous mapping,
which means that both images and inverse images of open sets are open
sets, and indeed the canonical quotient mapping is a nice covering
mapping, so that for every point $x$ in ${\bf R}^n$ there is a
neighborhood $U$ of $x$ in ${\bf R}^n$ such that the restriction of
$q$ to $U$ is a homeomorphism from $U$ onto the open set $q(U)$ in
${\bf R}^n / L$.  For that matter we can think of ${\bf R}^n / L$ as a
smooth manifold, with the quotient mapping $q$ as a smooth mapping
which is a local diffeomorphism.

	Suppose that $L_1$, $L_2$ are lattices in ${\bf R}^n$, and
let
\begin{equation}
	q_1 : {\bf R}^n \to {\bf R}^n / L_1, \quad
		q_2 : {\bf R}^n \to {\bf R}^n / L_2
\end{equation}
be the corresponding canonical quotient mappings.  If $A$ is an
invertible linear transformation on ${\bf R}^n$ such that
\begin{equation}
	A(L_1) = L_2,
\end{equation}
then we get an induced mapping
\begin{equation}
	\widehat{A} : {\bf R}^n / L_1 \to {\bf R}^n / L_2.
\end{equation}
This mapping is a group isomorphism and a homeomorphism, and even
a diffeomorphism, which satisfies the obvious compatibility condition
with the corresponding canonical quotient mappings $q_1$, $q_2$,
namely $q_1 \circ A = \widehat{A} \circ q_2$.

	When $n = 1$, one can consider the lattice $2 \pi {\bf Z}$
consisting of integer multiples of $2 \pi$, and it is customary to
identify ${\bf R} / 2 \pi {\bf Z}$ with the unit circle ${\bf T}$ in
the complex numbers ${\bf C}$,
\begin{equation}
	{\bf T} = \{z \in {\bf C} : |z| = 1\},
\end{equation}
where $|z|$ denotes the usual modulus of $z \in {\bf C}$, $|z| = (x^2
+ y^2)^{1/2}$ when $z = x + i \, y$, $x, y \in {\bf R}$.  More
precisely, $\exp (i \, t)$ is an explicit version of the canonical
quotient mapping from ${\bf R} / 2 \pi {\bf Z}$ onto ${\bf T}$ with
respect to this identification, which is a local diffeomorphism and a
group homomorphism using the group structure of multiplication on
${\bf T}$.  In general, we can identify ${\bf R}^n / 2 \pi {\bf Z}^n$
with ${\bf T}^n$, the $n$-fold Cartesian product of ${\bf T}$, where
$2 \pi {\bf Z}^n$ denotes the lattice of points whose coordinates are
all integer multiples of $2 \pi$.

	Suppose that $L$ is a lattice in ${\bf R}^n$.  Also let $A$ be
an invertible linear mapping on ${\bf R}^n$ such that $A(2 \pi {\bf
Z}^n) = L$.  Thus $\widehat{A}$ is a group isomorphism and a
diffeomorphism from ${\bf R}^n / 2 \pi {\bf Z}^n \cong {\bf T}^n$ onto
${\bf R}^n / L$.

	There is a more precise way to look at the quotient of ${\bf
R}^n$ by a lattice, which is to say that the quotient space has a kind
of local affine structure.  That is, there is a local affine structure
in which the canonical quotient mapping is considered to be locally
affine, and which permits one to say when a curve in the quotient is
locally a straight line segment, like an arc on a line, and when it
has locally constant speed, etc.  If $L_1$, $L_2$ are lattices in
${\bf R}^n$ and $A$ is an invertible linear mapping on ${\bf R}^n$
such that $A(L_1) = L_2$, then the induced mapping $\widehat{A}$ from
${\bf R}^n / L_1$ onto ${\bf R}^n / L_2$ preserves this local affine
structure on the quotient spaces.

	There is an even more precise way to look at the quotient
${\bf R}^n / L$ of ${\bf R}^n$ by a lattice $L$, which is that it has
a local flat geometric structure, induced from the one on ${\bf R}^n$.
With respect to this structure one can make local measurements of
lengths, volumes, and angles, like the length of a curve, the angle at
which two curves meet at a point, or the volume of a nice subset.  In
technical terms this can be seen as a \emph{Riemannian metric}.

	In particular, one can define the volume of such a quotient
${\bf R}^n / L$, where the volume of ${\bf R}^n / {\bf Z}^n$ is equal
to $1$, and the volume of ${\bf R}^n / 2 \pi {\bf Z}^n$ is equal to
$(2 \pi)^n$.  In general, if $L_1$, $L_2$ are lattices in ${\bf R}^n$
and $A$ is an invertible linear transformation on ${\bf R}^n$ such
that $A(L_1) = L_2$, then the volume of ${\bf R}^n / L_2$ is equal to
$|\det A|$ times the volume of ${\bf R}^n / L_1$, and more generally
if $E$ is a nice subset of ${\bf R}^n / L_1$, then the volume of
$\widehat{A}(E)$ in ${\bf R}^n / L_2$ is equal to $|\det A|$ times the
volume of $A$ in ${\bf R}^n / L_1$.  This is a variant of the fact
that on ${\bf R}^n$ a linear transformation $A$ distorts volumes by
a factor of $|\det A|$, where $\det A$ denotes the determinant of $A$.

	Suppose that $L_1$, $L_2$ are lattices in ${\bf R}^n$, and
that $T$ is an invertible linear transformation on ${\bf R}^n$ such
that $T(L_1) = L_2$.  Recall that $T$ is an \emph{orthogonal
transformation} on ${\bf R}^n$ if $T$ is invertible with inverse given
by the adjoint, also known as the transpose, of $T$, and that this is
equivalent to saying that $T$ preserves the standard norm of vectors
in ${\bf R}^n$, and the standard inner product of vectors in ${\bf
R}^n$.  In other words, orthogonal transformations on ${\bf R}^n$ are
linear mappings which preserve the geometry in ${\bf R}^n$, and for
the lattices $L_1$, $L_2$ and the quotients of ${\bf R}^n$ by them we
have that the induced mapping $\widehat{T}$ from ${\bf R}^n / L_1$
onto ${\bf R}^n / L_2$ preserves the geometry as well.

	In short, quotients of ${\bf R}^n$ by lattices are the same in
terms of group structure, topological and even smooth structure, and
affine structure, and not in general for more precise geometry.  The
volume of the quotient space is one basic parameter that one can
consider.  It is also interesting to look at closed curves in the
quotient which are locally flat, their lengths, the angles at which
they meet, and so on.

	We can consider lattices in ${\bf C}^n$ as well.  In this
regard we can identify ${\bf C}^n$ with ${\bf R}^{2n}$ in the usual
manner, so that the real and imaginary parts of the $n$ components of
an element of ${\bf C}^n$ give rise to the $2n$ components of an
element of ${\bf R}^{2n}$, and then define a lattice in ${\bf C}^n$ to
be a lattice in ${\bf R}^{2n} \sim {\bf C}^n$.  We write ${\bf Z}[i]$
for the \emph{Gaussian integers}, which are complex numbers of the
form $a + i \, b$, where $a$, $b$ are integers, and $({\bf Z}[i])^n$
for the lattice in ${\bf C}^n$ consisting of $n$-tuples of Gaussian
integers, which is also called the standard integer lattice in ${\bf
C}^n$.

	If $L$ is a lattice in ${\bf C}^n$, then the quotient ${\bf
C}^n / L$ inherits a complex structure from ${\bf C}^n$.  This means
in particular that the tangent spaces of the quotient are complex
vector spaces, just as they are for ${\bf C}^n$.  If $L_1$, $L_2$ are
lattices in ${\bf C}^n$ and $A$ is an invertible complex-linear
transformation on ${\bf C}^n$ such that $A(L_1) = L_2$, then $A$
induces a mapping $\widehat{A}$ from ${\bf C}^n / L_1$ to ${\bf C}^n /
L_2$ which preserves this complex structure.

	We can combine the complex and Riemannian structures and
consider \emph{Hermitian structures}.  Basically this means looking
at correspondences between lattices in ${\bf C}^n$ which come from
unitary mappings on ${\bf C}^n$.  If $L_1$, $L_2$ are lattices
in ${\bf C}^n$ and $T$ is a unitary mapping on ${\bf C}^n$
such that $T(L_1) = L_2$, then the induced mapping $\widehat{T}$
from ${\bf C}^n / L_1$ to ${\bf C}^n / L_2$ preserves both the
complex and Riemannian structures.

	Let us focus on complex structures for a moment.  It will be
convenient to write $\mathcal{L}({\bf R}^m, {\bf C}^n)$ for the space
of real-linear mappings from ${\bf R}^m$ to ${\bf C}^n$.  The complex
structure on ${\bf C}^n$ is still relevant for this space, in that
$\mathcal{L}({\bf R}^m, {\bf C}^n)$ is naturally a complex vector
space, because one can multiply elements of $\mathcal{L}({\bf R}^m,
{\bf C}^n)$ by $i$.  One can also describe these linear
transformations by $m \times n$ matrices of complex numbers in the
usual manner, using the standard bases for ${\bf R}^m$ and ${\bf
C}^n$.

	Let us write $\mathcal{L}^*({\bf R}^m, {\bf C}^n)$ for the
subset of $\mathcal{L}({\bf R}^m, {\bf C}^n)$ consisting of linear
transformations whose kernels are trivial, at least when $m \le 2n$,
so that this is possible.  Using the usual Euclidean topology for
$\mathcal{L}({\bf R}^m, {\bf C}^n)$, $\mathcal{L}^*({\bf R}^m, {\bf
C}^n)$ is an open set.  When $m = 2n$, $\mathcal{L}^*({\bf R}^m, {\bf
C}^n)$ consists of the invertible real-linear transformations from
${\bf R}^m$ onto ${\bf C}^n$, and a lattice in ${\bf C}^n$ is the
image of ${\bf Z}^{2n}$ under an element of $\mathcal{L}^*({\bf
R}^{2n}, {\bf C}^n)$.

	Now let us look at general lattices in ${\bf C}^n$, under the
equivalence relation in which two lattices $L_1$, $L_2$ are considered
to be equivalent if there is an invertible complex-linear
transformation $A$ on ${\bf C}^n$ such that $A(L_1) = L_2$.  This
leads to an equivalence relation on $\mathcal{L}^*({\bf R}^{2n}, {\bf
C}^n)$, in which two elements of $\mathcal{L}^*({\bf R}^{2n}, {\bf
C}^n)$ are considered to be equivalent if one can be written as the
composition of an invertible complex-linear transformation on ${\bf
C}^n$ with the other element of $\mathcal{L}^*({\bf R}^{2n}, {\bf
C}^n)$.  In other words, we look at the action of $GL({\bf C}^n)$ on
$\mathcal{L}^*({\bf R}^{2n}, {\bf C}^n)$ by post-composition.

	Actually, it is more convenient to consider
$\mathcal{L}_1^*({\bf R}^{2n}, {\bf C}^n)$, which we define to be the
subset of $\mathcal{L}^*({\bf R}^{2n}, {\bf C}^n)$ consisting of
invertible real-linear transformations from ${\bf R}^{2n}$ to ${\bf
C}^n$ such that the image of the first $n$ standard basis vectors in
${\bf R}^{2n}$ are linearly-independent over the complex numbers as
$n$ vectors in ${\bf C}^n$.  This restriction is not too serious, and
indeed we can describe the lattices in ${\bf C}^n$ as images of ${\bf
Z}^{2n}$ under mappings in $\mathcal{L}_1^*({\bf R}^{2n}, {\bf C}^n)$.
In other words, if we start with a lattice $L$ given as the image of
${\bf Z}^{2n}$ under an element of $\mathcal{L}^*({\bf R}^n, {\bf
C}^n)$, we can rewrite it as the image of ${\bf Z}^{2n}$ under a
linear transformation in $\mathcal{L}_1^*({\bf R}^{2n}, {\bf C}^n)$ by
pre-composing the initial linear transformation from ${\bf R}^{2n}$ to
${\bf C}^n$ with an invertible linear transformation on ${\bf R}^{2n}$
which permutes the standard basis vectors in a suitable way.

	To deal with the action of $GL({\bf C}^n)$ by
post-composition, we can restrict ourselves to $\mathcal{L}^{**}({\bf
R}^{2n}, {\bf C}^n)$, which we define to be the space of invertible
real-linear transformations from ${\bf R}^{2n}$ to ${\bf C}^n$ such
that the images of the first $n$ standard basis vectors in ${\bf
R}^{2n}$ are the $n$ standard basis vectors in ${\bf C}^n$, and in the
same order.  In other words, if we identify ${\bf R}^{2n}$ with the
Cartesian product ${\bf R}^n \times {\bf R}^n$, then these are the
invertible real-linear transformations from ${\bf R}^n \times {\bf
R}^n$ onto ${\bf C}^n$ with the property that on ${\bf R}^n \times
\{0\}$ they coincide with the standard embedding of ${\bf R}^n$ into
${\bf C}^n$.  This exactly compensates for the action of $GL({\bf C}^n)$
by post-composition, since for any collection $v_1, \ldots, v_n$
of linearly-independent vectors in ${\bf C}^n$ there is a unique
$A \in GL({\bf C}^n)$ such that $A(v_1), \ldots, A(v_n)$
are the standard basis vectors in ${\bf C}^n$, in order.

	We can identify $\mathcal{L}^{**}({\bf R}^{2n}, {\bf C}^n)$
with an open subset of $\mathcal{L}({\bf R}^n, {\bf C}^n)$.  That is,
elements of $\mathcal{L}^{**}({\bf R}^{2n}, {\bf C}^n)$ can be
identified with linear transformations from ${\bf R}^n \times {\bf
R}^n$ into ${\bf C}^n$, and these linear transformations are
determined by what they do on $\{0\} \times {\bf R}^n$, since their
behavior on ${\bf R}^n \times \{0\}$ is fixed by definition.  We can
think of elements of $\mathcal{L}({\bf R}^n, {\bf C}^n)$ as being
written as $A + i \, B$, where $A$, $B$ are linear transformations on
${\bf R}^n$, and one can check that the elements of
$\mathcal{L}^{**}({\bf R}^{2n}, {\bf C}^n)$ correspond exactly to
elements of $\mathcal{L}({\bf R}^n, {\bf C}^n)$ of the form $A + i \,
B$, where $A$, $B$ are linear transformations on ${\bf R}^n$ and $B$
is invertible.

	To be more precise, it is helpful to think in terms of 
real-linear mappings on ${\bf C}^n$, which can be written as
\begin{equation}
	T(x + i \, y) = E_1(x) + E_2(y) + i (E_3(x) + E_4(y)),
\end{equation}
where $x, y \in {\bf R}^n$.  The passage to $\mathcal{L}_1^*({\bf
R}^{2n}, {\bf C}^n)$ can be expressed in these terms as the
restriction to invertible real-linear transformations $T$ on ${\bf
C}^n$ of the form
\begin{equation}
	T(x + i \, y) = x + A(y) + i \, B(y),
\end{equation}
where $A$, $B$ are linear transformations on ${\bf R}^n$.  The
condition of invertibility of $T$ is equivalent to the invertibility
of $B$ on ${\bf C}^n$.

	Another way to look at real-linear mappings on ${\bf C}^n$
is as mappings of the form
\begin{equation}
	T(z) = M(z) + \overline{N(z)},
\end{equation}
where $z \in {\bf C}^n$, $M$ and $N$ are complex-linear mappings on
${\bf C}^n$, and for $w \in {\bf C}^n$, $\overline{w}$ is the element
of ${\bf C}^n$ whose coordinates are the complex-conjugates of the
coordinates of $w$.

	Invertibility of $T$ is a bit tricky, and as an
important special case, it is natural to restrict our attention to
mappings $T$ as above for which $M$ majorizes $N$ in the sense that
\begin{equation}
	|N(z)| < |M(z)|
\end{equation}
for $z \in {\bf C}^n$, $z \ne 0$, where $|w|$ denotes the standard
Euclidean norm of $w \in {\bf C}^n$.  To factor out the action of
$GL({\bf C}^n)$ by post-composition, we can restrict our attention to
real-linear transformations $T$ of the form
\begin{equation}
	T(z) = z + \overline{E(z)},
\end{equation}
where $E$ is a complex-linear transformation on ${\bf C}^n$ with
operator norm strictly less than $1$, which is equivalent to saying
that $E^* \, E < I$.  This has nice features when we think of the
image of the standard integer lattice $({\bf Z}[i])^n$ under $T$, with
points in the image being reasonably-close to their counterparts in
the original lattice.

	The $n = 1$ case is quite instructive.  We can write a
real-linear transformation $T$ on ${\bf C}$ as
\begin{equation}
	T(x + i \, y) = a \, x + i \, b \, y
\end{equation}
for $x, y \in {\bf R}$, where $a$, $b$ are complex numbers, and when
$T$ is invertible we can rewrite this as
\begin{equation}
	T(x + i \, y) = a (x + i \, c \, y),
\end{equation}
where $a$, $c$ are complex numbers with $a \ne 0$ and $c$ having
nonzero imaginary part.  Alternatively, we can write a real-linear
transformation $T$ on ${\bf C}$ as $T(z) = \alpha \, z + \beta \,
\overline{z}$ with $\alpha, \beta \in {\bf C}$, and where $T$ is
invertible if and only if $|\alpha| \ne |\beta|$, and when
$|\alpha| > |\beta|$ this can be rewritten as
\begin{equation}
	T(z) = \theta (z + \mu \, \overline{z}),
\end{equation}
where $\theta$ is a nonzero complex number and $\mu$ is a complex
number such that $|\mu| < 1$.

	Let us return now to the real case.  Consider the quotient
space $O({\bf R}^n) \backslash GL({\bf R}^n)$, in which two invertible
linear transformations on ${\bf R}^n$ are identified if one can be
written as an orthogonal linear transformation times the other.  We
can identify this quotient space with the space of symmetric linear
transformations on ${\bf R}^n$ which are positive definite, through
the mapping
\begin{equation}
	T \mapsto T^* \, T.
\end{equation}
In other words, if $T$ is an invertible linear transformation on ${\bf
R}^n$, then $T^* \, T$ is a symmetric linear transformation on ${\bf
R}^n$ which is positive-definite, $T_1^* \, T_1 = T_2^* \, T_2$ for
$T_1, T_2 \in GL({\bf R}^n)$ if and only if $T_2 = R \, T_1$ for some
orthogonal transformation $R$, and every symmetric linear transformation
on ${\bf R}^n$ which is positive-definite can be expressed as $T^* \, T$
for an invertible linear transformation $T$.

	Similarly, the quotient $SO({\bf R}^n) \backslash SL({\bf
R}^n)$ can be identified with the space $\mathcal{M}({\bf R}^n)$ of
symmetric linear transformations on ${\bf R}^n$ which are positive
definite and have determinant equal to $1$.  Let us write $\Sigma({\bf
R}^n)$ for the elements of $SL({\bf R}^n)$ whose matrices with respect
to the standard basis have integer entries.  The inverse of a linear
transformation in $\Sigma({\bf R}^n)$ also lies in $\Sigma({\bf
R}^n)$, because Cramer's rule gives a formula for the matrix of the
inverse which shows that it has integer entries when the original
matrix has integer entries and determinant equal to $1$.  

	Elements of $\Sigma({\bf R}^n)$ can be described as the
invertible linear transformations which take ${\bf Z}^n$ onto itself.
The quotient $SL({\bf R}^n) / \Sigma({\bf R}^n)$ describes the space
of lattices $L$ in ${\bf R}^n$ such that the corresponding quotient
${\bf R}^n / L$ has volume equal to $1$ and for which there is an
extra piece of data concerning orientation, and the double quotient
$SO({\bf R}^n) \backslash SL({\bf R}^n) / \Sigma({\bf R}^n)$ deals
with these lattices up to equivalence under rotation.  By identifying
$SO({\bf R}^n) \backslash SL({\bf R}^n)$ with $\mathcal{M}({\bf
R}^n)$, the double quotient can be identified with the quotient of
$\mathcal{M}({\bf R}^n)$ by the action of $\Sigma({\bf R}^n)$ defined
by $A \mapsto T^* \, A \, T$, $A \in \mathcal{M}({\bf R}^n)$, $T \in
\Sigma({\bf R}^n)$.  This quotient is denoted $\mathcal{M}({\bf R}^n)
/ \Sigma({\bf R}^n)$.

	In the complex case let us consider lattices $L$ in ${\bf
C}^n$ which are of the form $A(({\bf Z}[i])^n)$ for some invertible
complex-linear mapping $A$ on ${\bf C}^n$.  It is natural to look at
these lattices up to unitary equivalence, which is to say that two
lattices $L_1$, $L_2$ are equivalent if there is a unitary linear
transformation $T$ on ${\bf C}^n$ such that $T(L_1) = L_2$.  This
leads to an equivalence relation on $GL({\bf C}^n)$, in which two
invertible linear transformations $A_1$, $A_2$ on ${\bf C}^n$ are
considered to be equivalent if there is a unitary linear
transformation $T$ on ${\bf C}^n$ such that $A_2 = T \, A_1$.  The
quotient of $GL({\bf C}^n)$ by this equivalence relation is denoted
\begin{equation}
	U({\bf C}^n) \backslash GL({\bf C}^n)
\end{equation}
and can be identified with the space of self-adjoint linear
transformations on ${\bf C}^n$ which are positive definite, through
the mapping
\begin{equation}
	A \in GL({\bf C}^n) \mapsto A^* \, A.
\end{equation}
That is, for each element $A$ of $GL({\bf C}^n)$, the product $A^* \,
A$ is a self-adjoint linear transformation on ${\bf C}^n$ which is
positive-definite, $A_1^* \, A_1 = A_2^* \, A_2$ for two elements
$A_1$, $A_2$ of $GL({\bf C}^n)$ if and only if $A_2 = T \, A_1$ for
some unitary linear transformation $T$ on ${\bf C}^n$, and every
self-adjoint linear transformation on ${\bf C}^n$ can be expressed as
$A^* \, A$ for some invertible linear transformation $A$ on ${\bf
C}^n$.

	Similarly, one can consider two elements $B_1$, $B_2$ of
$SL({\bf C}^n)$ to be equivalent when there is a linear transformation
$U$ in the special unitary group $SU({\bf C}^n)$ such that $A_2 = U \,
A_1$.  The quotient $SU({\bf C}^n) \backslash SL({\bf C}^n)$ can be
identified with the space of self-adjoint linear transformations on
${\bf C}^n$ which are positive-definite and have determinant $1$,
through the same mapping as before.  Let us consider lattices $L$ of
the form $B(({\bf Z}[i])^n)$ for some $B \in SL({\bf C}^n)$, a modest
normalization.

	As in the real case we write $\Sigma({\bf C}^n)$ for the
subgroup of $SL({\bf C}^n)$ of linear transformations whose associated
$n \times n$ matrices, with respect to the standard basis for ${\bf
C}^n$, have integer entries, which implies that the matrices
associated to their inverses also have integer entries.  Thus $B(({\bf
Z}[i])^n) = ({\bf Z}[i])^n$ when $B \in \Sigma({\bf R}^n)$, and
conversely $B \in SL({\bf C}^n)$ and $B(({\bf Z}[i])^n) = ({\bf
Z}[i])^n$ implies that $B \in \Sigma({\bf C}^n)$.  The quotient
$SL({\bf C}^n) / \Sigma({\bf C}^n)$ represents the space of lattices
under consideration, the double quotient $SU({\bf C}^n) \backslash
SL({\bf C}^n) / \Sigma({\bf C}^n)$ represents the space of these
lattices modulo equivalence under special unitary transformations, and
this double quotient can be identified with the quotient of the space
$\mathcal{M}({\bf C}^n)$ of self-adjoint positive-definite linear
transformations on ${\bf C}^n$ with determinant $1$ by the action of
$\Sigma({\bf C}^n)$ defined by $P \mapsto B^* \, P \, B$, $B \in
\Sigma({\bf C}^n)$.

	Next we consider real and complex projective spaces.  Namely,
if $n$ is a positive integer, then the $n$-dimensional real and
complex projective spaces ${\bf RP}^n$, ${\bf CP}^n$ consist of the
real and complex lines through the origin in ${\bf R}^{n+1}$, ${\bf
C}^{n+1}$, respectively.  To put it another way, if ${\bf R}_*$, ${\bf
C}_*$ denote the nonzero real and complex numbers, respectively, then
we have natural actions of ${\bf R}_*$, ${\bf C}_*$ on ${\bf R}^{n+1}
\backslash \{0\}$, ${\bf C}^{n+1} \backslash \{0\}$ by scalar
multiplication, and the projective spaces are the corresponding
quotient spaces.  Thus two nonzero vectors $v$, $w$ in ${\bf
R}^{n+1}$, ${\bf C}^{n+1}$ lead to the same point in the corresponding
projective space exactly when they are scalar multiples of each other.
Note that we get canonical mappings from ${\bf R}^{n+1} \backslash
\{0\}$, ${\bf C}^{n+1} \backslash \{0\}$ onto ${\bf RP}^n$, ${\bf
CP}^n$, in which a nonzero vector $v$ is sent to the line through the
origin which passes through $v$.

	If $L$ is a nontrivial linear subspace of ${\bf R}^{n+1}$,
${\bf C}^{n+1}$ of dimension $l + 1$, say, then we get an interesting
space ${\bf P}(L)$ consisting of all lines through the origin in $L$,
which we can think of as sitting inside of ${\bf RP}^n$, ${\bf CP}^n$,
as appropriate.  More precisely, ${\bf P}(L)$ is basically a copy of
${\bf RP}^l$ or ${\bf CP}^l$.  These are the $l$-dimensional ``linear
subspaces'' of projective space, analogous to linear subspaces of
${\bf R}^n$, ${\bf C}^n$.

	If $A$ is an invertible linear transformation on ${\bf
R}^{n+1}$ or on ${\bf C}^{n+1}$, then $A$ takes lines to lines, and
induces a transformation $\widehat{A}$ on the corresponding projective
space.  Notice that $\widehat{A}$ is automatically a one-to-one
transformation of the corresponding projective space onto itself, with
\begin{equation}
	\widehat{A}^{-1} = \widehat{(A^{-1})},
\end{equation}
and $\widehat{A}$ maps linear subspaces of projective space to
themselves, in the sense of the preceding paragraph.  Also, if $A_1$,
$A_2$ are invertible linear transformations on ${\bf R}^{n+1}$ or on
${\bf C}^{n+1}$, then the induced transformations $\widehat{A}_1$,
$\widehat{A}_2$ on the corresponding projective space satisfy
\begin{equation}
	\widehat{A}_1 \circ \widehat{A}_2 = \widehat{(A_1 \circ A_2)}.
\end{equation}

	Let $H$ be a hyperplane in ${\bf R}^{n+1}$ or in ${\bf
C}^{n+1}$, which is to say a linear subspace of dimension $n$, and let
$v$ be a nonzero vector in ${\bf R}^{n+1}$, ${\bf C}^{n+1}$, as
appropriate.  This leads to an affine hyperplane $H + v$, consisting
of all vectors of the form $w + v$, $w \in H$, and which does not
contain the vector $0$.  For each $w \in H$, we can look at the line
through $w + v$, which we can view as an element of the corresponding
projective space.

	In other words, we basically get an embedding of $H$ into the
corresponding projective space, ${\bf RP}^n$ or ${\bf CP}^n$.  Of
course we can also think of $H$ as being isomorphic to ${\bf R}^n$ or
${\bf C}^n$, so that we are really looking at a bunch of embeddings of
${\bf R}^n$, ${\bf C}^n$ into ${\bf RP}^n$, ${\bf CP}^n$,
respectively.  For instance, we can do this with $H$ equal to the
$j$th coordinate hyperplane in ${\bf R}^{n+1}$, ${\bf C}^{n+1}$,
$1 \le j \le n + 1$, which is defined by the condition that the
$j$th coordinate of vectors in $H$ are equal to $0$, and we can take
$v$ to be the $j$th standard basis vector, with $j$th coordinate
equal to $1$ and the other $n$ coordinates equal to $0$.

	These $n + 1$ embeddings of ${\bf R}^n$, ${\bf C}^n$ into
${\bf RP}^n$, ${\bf CP}^n$ corresponding to the $n + 1$ coordinate
hyperplanes in ${\bf R}^{n+1}$, ${\bf CP}^{n+1}$ are sufficient to
cover the projective space, i.e., every point in projective space
shows up in the image of at least one of the embeddings.  For a given
hyperplane $H$, the set of points in the projective space which do not
occur in the embedding of $H$ is the same as ${\bf P}(H)$.  Thus
the set of missing points in the projective space lie in a projective
subspace of dimension $1$ less.

	Using these embeddings of ${\bf R}^n$, ${\bf C}^n$ into the
corresponding projective spaces, we can think of the projective spaces
as being manifolds.  That is, these embeddings provide local
coordinates for all points in the projective space.  Two different
embeddings which contain the same point $p$ in their image are
compatible in terms of topology and also smooth structure, in the
complex case we can say that ${\bf CP}^n$ is a complex manifold.  In
the real and complex situations there is a finer ``projective''
structure which is reflected in the presence of nice projective
subspaces, for instance, and the projectivized versions of linear
transformations on ${\bf R}^{n+1}$, ${\bf C}^{n+1}$.

	Note that two invertible linear transformations $A_1$, $A_2$
on ${\bf R}^{n+1}$ or on ${\bf C}^{n+1}$ lead to the same induced
transformation on projective space if and only if there is a nonzero
scalar $\alpha$ such that $A_2 = \alpha \, A_1$.  Thus the group of
these ``projective linear transformations'' has dimension $(n+1)^2 -
1$ over the real or complex numbers, as appropriate.  Also, for
any pair of points $p$, $q$ in a projective space, there is a
projective linear transformation which takes $p$ to $q$.

	If $k$, $n$ are positive integers with $k < n$, then the
Grassmann spaces $G_{\bf R}(k, n)$, $G_{\bf C}(k, n)$ consist of the
$k$-dimensional linear subspaces of ${\bf R}^n$, ${\bf C}^n$,
respectively.  When $k = 1$ this reduces to the $(n-1)$-dimensional
projective spaces.  Suppose that $L$, $M$ are linear subspaces of
${\bf R}^n$ or of ${\bf C}^n$ which are complementary, in the sense
that
\begin{equation}
	L \cap M = \{0\} \quad\hbox{and}\quad 
		L + M = {\bf R}^n \hbox{ or } {\bf C}^n,
\end{equation}
as appropriate, and that $L$ has dimension $k$, so that $M$ has
dimension $n - k$.  If $A$ is a linear mapping from $L$ to $M$, then
the graph of $A$, consisting of the vectors
\begin{equation}
	v + A(v), \quad v \in L,
\end{equation}
is also a $k$-dimensional subspace of ${\bf R}^n$ or of ${\bf C}^n$,
as appropriate.  In this way we can embed the vector space of linear
transformations from $L$ to $M$ into the Grassmannian, and this
provides a nice coordinate patch around $L$ itself.

	In particular, these coordinate patches permit one to view the
Grassmann spaces as smooth manifolds, and as complex manifolds in the
complex case.  The dimension of $G_{\bf R}(k, n)$, $G_{\bf C}(k, n)$
is equal to
\begin{equation}
	k (n-k)
\end{equation}
with respect to the real or complex numbers, as appropriate.  Just as
for projective spaces, invertible linear transformations on ${\bf
R}^n$ or on ${\bf C}^n$ induce interesting mappings on the
corresponding Grassmannians.  These actions are again transitive,
because if $L_1$, $L_2$ are $k$-dimensional linear subspaces of ${\bf
R}^n$ or of ${\bf C}^n$, then there is an invertible linear
transformation $A$ on ${\bf R}^n$ or on ${\bf C}^n$, as appropriate,
such that $A(L_1) = L_2$.

	There is a natural correspondence between the Grassmann spaces
of $k$-dimensional linear subspaces in ${\bf R}^n$ or in ${\bf C}^n$
and the Grassmann spaces of $(n-k)$-dimensional linear subspaces of
${\bf R}^n$ or in ${\bf C}^n$, respectively.  More precisely, there is
a natural correspondence between $k$-dimensional linear subspaces of
${\bf R}^n$ or of ${\bf C}^n$ and $(n - k)$-dimensional linear
subspaces of the dual spaces associated to ${\bf R}^n$ or ${\bf C}^n$,
consisting of the linear functionals on ${\bf R}^n$ or on ${\bf C}^n$.
Namely, if $L$ is a $k$-dimensional linear subspace of ${\bf R}^n$ or
of ${\bf C}^n$, then one gets an $(n-k)$-dimensional linear subspace
of the dual space by taking the linear functionals which vanish on
$L$.  Conversely, if one starts with an $(n-k)$-dimensional subspace
of the dual space, one gets a $k$-dimensional subspace of the original
space by taking the intersections of the kernels of the linear
functionals in the subspace of the dual space.  Using the standard
basis for ${\bf R}^n$ or ${\bf C}^n$, or any other basis for that matter,
one can identify ${\bf R}^n$, ${\bf C}^n$ with their dual spaces
in a well-known manner.

	Suppose that $L$ is a linear subspace of ${\bf R}^n$ or of
${\bf C}^n$ of dimension $l$.  If $l \ge k$, then we get an
interesting subspace of the Grassmann space of $k$-dimensional linear
subspaces of ${\bf R}^n$ or ${\bf C}^n$, as appropriate, consisting of
the $k$-dimensional linear subspaces contained in $L$.  If $k \ge l$,
then there is another interesting subspace of the Grassmann space,
consisting of the $k$-dimensional linear subspaces of ${\bf R}^n$ or
${\bf C}^n$ which contain $L$.  When $k = l$, these two cases are the
same and we get simply a point in the Grassmann space.

	Let $n$ be a positive integer.  By a \emph{multiindex}
we mean an $n$-tuple $\alpha = (\alpha_1, \ldots, \alpha_n)$
of nonnegative integers, and in this case we set
\begin{equation}
	|\alpha| = \sum_{j=1}^n |\alpha_j|.
\end{equation}
Given a multiindex $\alpha$, we define the corresponding
monomial $w^\alpha$ on ${\bf R}^n$ or on ${\bf C}^n$
by
\begin{equation}
	w^\alpha = w_1^{\alpha_1} \cdots w_n^{\alpha_n},
\end{equation}
where $w = (w_1, \ldots, w_n)$ as usual, and we call $|\alpha|$ the
\emph{degree} of this monomial.  When $\alpha_j = 0$ we interpret
$w_j^{\alpha_j}$ as being equal to the constant $1$.

	A polynomial $p(w)$ on ${\bf R}^n$ or on ${\bf C}^n$ is a
function which is a linear combination of monomials.  We take the
coefficients to be real or complex numbers according to whether we are
working on ${\bf R}^n$ or on ${\bf C}^n$.  A polynomial $p(w)$ which
is a linear combination of monomials of the same degree $a$ is said to
be homogeneous of degree $a$.  This is equivalent to the condition
that
\begin{equation}
	p(\lambda \, w) = \lambda^a \, p(w)
\end{equation}
for all real or complex numbers $\lambda$ and all $w$ in ${\bf R}^n$
or in ${\bf C}^n$, as appropriate.

	Suppose that $p_1(w), \ldots, p_{n+1}(w)$ are polynomials on
${\bf R}^{n+1}$ or on ${\bf C}^{n+1}$ which are homogeneous of the
same degree $a > 0$.  Assume also that
\begin{equation}
	p_1(w) = \cdots = p_{n+1}(w) = 0
\end{equation}
only when $w = 0$.  The combined mapping $p(w) = (p_1(w), \ldots,
p_{n+1}(w))$ is a homogeneous polynomial mapping of ${\bf R}^{n+1}$ or
${\bf C}^{n+1}$ into itself which maps nonzero vectors to nonzero
vectors, and as a result induces a mapping $\widehat{p}$ from ${\bf
RP}^n$ or ${\bf CP}^n$ to itself, as appropriate.  The degree $1$ case
corresponds exactly to invertible linear mappings on ${\bf R}^{n+1}$
or on ${\bf C}^{n+1}$ and the associated projective linear
transformations on the corresponding projective spaces.  If $p = (p_1,
\ldots, p_{n+1})$, $q = (q_1, \ldots, q_{n+1})$ are homogeneous
polynomial mappings on ${\bf R}^{n+1}$ or on ${\bf C}^{n+1}$ of this
type, of degrees $a, b > 0$, respectively, then the composition $p
\circ q$ is a homogeneous polynomial mapping of degree $a b$, and
$\widehat{p \circ q} = \widehat{p} \circ \widehat{q}$.

	When $n = 1$, we can think of ${\bf RP}^1$, ${\bf CP}^1$ as
being the same as ${\bf R}$, ${\bf C}$ with an additional point added,
often denoted $\infty$.  Projective linear transformations on ${\bf
RP}^1$, ${\bf CP}^1$ can then be described as mappings of the
form $(a \, z + b)/(c \, z + d)$, where $a \, d - b \, c \ne 0$,
and using standard conventions along the lines of $1 / 0 = \infty$,
$1 / \infty = 0$.  Similarly, the mappings associated to homogeneous
polynomials as in the preceding paragraph reduce to nonconstant
rational functions of a single variable.

\section{Immersions, submersions, and connections}
\label{Immersions, submersions, and connections}
\setcounter{equation}{0}

	Let $M$, $N$ be nonempty $m$, $n$-dimensional smooth
manifolds, so that $M$, $N$ look locally like ${\bf R}^m$, ${\bf R}^n$
in a sense, and they have countable bases for their topologies.  As
basic situations, $M$, $N$ might in fact be open subsets of ${\bf
R}^m$, ${\bf R}^n$, respectively.  They might instead be given as
embedded submanifolds of higher-dimensional Euclidean spaces.

	Suppose that $f$ is a smooth mapping from $M$ into $N$.  For
each point $p \in M$, the differential of $f$ at $p$ is denoted $df_p$
and is a linear mapping from the tangent space of $M$ at $p$, which is
denoted $T_p M$ and is an $m$-dimensional real vector space, into the
tangent space $T_{f(p)} N$ of $N$ at $f(p)$, an $n$-dimensional real
vector space.  We say that $f$ is an \emph{immersion} if $df_p$ is an
injective linear transformation from $T_p M$ into $T_{f(p)} N$ for
each $p \in M$, which is to say that the kernel of $df_p$ is trivial
for all $p \in M$, and in which case $m \le n$.  We say that $f$ is a
\emph{submersion} if $df_p$ maps $T_p M$ onto $T_{f(p)} N$ for all $p
\in M$, in which case $m \ge n$.  When $m = n$, these two conditions
are equivalent to each other, and to the statement that $df_p$ is a
one-to-one linear mapping from $T_p M$ onto $T_{f(p)} N$ for each $p
\in M$.

	When $m \le n$, we have the standard embedding of ${\bf R}^m$
into ${\bf R}^n$, in which a point $x = (x_1, \ldots, x_m)$ in ${\bf
R}^m$ is sent to $\widehat{x} = (\widehat{x}_1, \ldots,
\widehat{x}_n)$ in ${\bf R}^n$, with $\widehat{x}_i = x_i$ for $1 \le
i \le m$ and $\widehat{x}_i = 0$ for $m < i \le n$.  When $m \ge n$,
we have the standard projection from ${\bf R}^m$ onto ${\bf R}^n$, in
which one keeps the first $n$ coordinates of a point in ${\bf R}^m$
and drops the remaining $m - n$ coordinates.  By the implicit function
theorem, immersions and submersions are locally equivalent to these
standard models, with the appropriate dimensions.  When $m = n$,
immersions and submersions are the same, and they are local
diffeomorphisms, as in the inverse function theorem.

	If $f : M \to N$ is a submersion, then $f$ is an open mapping
in particular, which is to say that $f(U)$ is an open subset of $N$
whenever $U$ is an open subset of $M$.  If we also assume that $f$ is
proper, in the sense that $f^{-1}(K)$ is a compact subset of $M$
whenever $K$ is a compact subset of $N$, then it is easy to check that
$f(E)$ is a closed subset of $N$ when $E$ is a closed subset of $M$.
This implies in turn that $f(M) = N$ when $N$ is connected, since
$f(M)$ would then be a nonempty subset of $N$ which is both open and
closed.

	Let us consider some examples.  Fix a positive integer $n$,
and for $M$ take ${\bf R}^{n+1} \backslash \{0\}$.  For $N$ we can
take the projective space ${\bf RP}^n$, and we have a natural mapping
from ${\bf R}^{n+1} \backslash \{0\}$ to ${\bf RP}^n$ which sends a
nonzero vector $v$ in ${\bf R}^{n+1}$ to the point in ${\bf RP}^n$
corresponding to the line through $v$.  This mapping is smooth and
defines a submersion.

	We can also restrict this mapping to the unit sphere ${\bf
S}^n$ in ${\bf R}^{n+1}$, consisting of the vectors $v$ such that $|v|
= 1$.  This mapping is still a smooth mapping from ${\bf S}^n$ onto
${\bf RP}^n$, with the dimensions of the domain and range now being
the same.  This mapping is a local diffeomorphism, and of course two
vectors $v$, $w$ in ${\bf S}^n$ are sent to the same point in ${\bf
RP}^n$ if and only if $v = w$ or $v = -w$.

	Now take $M$ to be ${\bf C}^{n+1} \backslash \{0\}$, which has
real dimension $2n + 2$, and let $N$ be ${\bf CP}^n$, which has real
dimension $2n$.  Once again there is a natural mapping which sends a
nonzero vector $v$ in ${\bf C}^{n+1}$ to the point in complex
projective space ${\bf CP}^n$ that corresponds to the line through
$v$.  This mapping is smooth, and in fact holomorphic, and it is also
a submersion.  We can also restrict this mapping to the sphere ${\bf
S}^{2n+1}$ consisting of the vectors $v$ in ${\bf C}^{n+1}$ such that
$|v| = 1$, to get a submersion onto ${\bf CP}^n$.  Two vectors $v$, $w$
in ${\bf S}^{2n+1}$ are sent to the same point in ${\bf CP}^n$ by this
mapping if and only if $w = \alpha \, v$ for some complex number $\alpha$
such that $|\alpha| = 1$, so that the fibers of this submersion from
${\bf S}^{2n+1}$ onto ${\bf CP}^n$ are all circles.

	Suppose that $M$, $N$ are smooth manifolds of dimensions $m$,
$n$, respectively, and that $f : M \to N$ is a proper smooth
submersion, so that the fibers $f^{-1}(z)$, $z \in N$, are compact
submanifolds of $M$ of dimension $m - n$.  Fix a point $z_1$ in $N$,
and suppose that $V_1$ is a neighborhood of $z_1$ in $N$ and that
$\phi$ is a smooth mapping from $f^{-1}(V_1) \subseteq M$ into
$f^{-1}(z_1)$ such that $\phi(x) = x$ when $f(x) = z_1$.  We can
combine $f$, $\phi$ to get a smooth mapping $(f, \phi)$ from
$f^{-1}(V_1)$ into $V_1 \times f^{-1}(z_1)$.  The differential of this
combined mapping is invertible at each element of the fiber
$f^{-1}(z_1)$, and it follows that there is a neighborhood $V_2$ of
$z_1$ contained in $V_1$ such that the combined mapping $(f, \phi)$
defines a diffeomorphism from $f^{-1}(V_2)$ onto $V_2 \times
f^{-1}(z_1)$.  In particular, if $z_2$ is an element of $N$ which is
sufficiently close to $z_1$, then the fibers $f^{-1}(z_1)$, $f^{-1}(z_2)$
are diffeomorphic smooth manifolds of dimension $m - n$.  If $N$ is
connected, then it follows that all of the fibers $f^{-1}(z)$, $z \in N$,
are diffeomorphic to each other.

	Again let $M$, $N$ be smooth manifolds of dimensions $m$, $n$,
respectively, and let $f : M \to N$ be a smooth submersion which may
or may not be proper, at least for the moment.  For each $p \in M$, we
get a linear subspace $V_p$ of the tangent space $T_p M$ consisting of
the tangent vectors to $M$ at $p$ which are also tangent to the fiber
$f^{-1}(f(p))$ of $f$ passing through $p$.  This can also be described
as the kernel of the differential $df_p$ of $f$ at $p$, as a linear
mapping from $T_p M$ to $T_{f(p)} N$.  Of course $V_p$ has dimension
$m - n$ for each $p \in M$, because $f$ is a submersion.

	By a \emph{connection} on $M$ with respect to the submersion
$f : M \to N$ we mean a choice of an $n$-dimensional linear subspace
$H_p$ of the tangent space $T_p M$ for each point $p \in M$ which
is transverse to the vertical subspace $V_p$ and which depends
smoothly on $p$.  We call $H_p$ the horizontal linear subspace of
the tangent space $T_p M$ determined by the connection, and the
restriction of the differential $df_p$ of $f$ at $p$ to the horizontal
subspace $H_p$ of $T_p M$ is a one-to-one linear mapping of $H_p$
onto the tangent space $T_{f(p)} N$ of $N$ at $f(p)$.  To put it another
way, at each point $p \in M$ the tangent space $T_p M$ of $M$ at $p$
is the direct sum of the horizontal and vertical subspaces $H_p$ and
$V_p$.  The vertical subspace $V_p$ is determined by $f$, and there is
some room for making choices for the horizontal subspaces.

	As a basic scenario, suppose that $M$ is also equipped with a
smooth Riemannian metric, which is to say an inner product on each
tangent space $T_p M$ which depends smoothly on $p$.  In this case,
one can simply choose $H_p$ to be the orthogonal complement of $V_p$
in $T_p M$ with respect to the Riemannian metric.  This shows that
connections always exist, since Riemannian metrics exist on any smooth
manifold.  Recall that one way to choose a Riemannian metric on a smooth
manifold $M$ is to choose local Riemannian metrics in coordinate charts
and combine them using a partition of unity, and another way is to embed
$M$ into a Euclidean space and then use the Riemannian metric inherited
from the one on the Euclidean space.

	Let $A$ be an invertible linear transformation on ${\bf R}^n$
or on ${\bf C}^n$ such that
\begin{equation}
	\lim_{l \to \infty} A^l(v) = 0
\end{equation}
for all $v$ in ${\bf R}^n$ or ${\bf C}^n$, as appropriate.  A
sufficient condition for this to hold is that the norm of $A$ be
strictly less than $1$.  In the complex case, this condition holds if
and only if the eigenvalues of $A$ all have modulus strictly less than
$1$, and in the real case one can complexify ${\bf R}^n$ to convert
$A$ to a linear transformation on ${\bf C}^n$, and the condition holds
on ${\bf R}^n$ if and only if it holds for the complexification on
${\bf C}^n$, which is to say that the modulus of each of the
eigenvalues of the associated linear transformation on ${\bf C}^n$
should be strictly less than $1$.  Notice that our condition is also
equivalent to
\begin{equation}
	\lim_{l \to \infty} |A^{-l}(v)| = \infty
\end{equation}
for all nonzero vectors $v$ in ${\bf R}^n$ or ${\bf C}^n$, as
appropriate.

	Let us define a space $\mathcal{H}_A$ by starting with ${\bf
R}^n \backslash \{0\}$ or ${\bf C}^n \backslash \{0\}$, as
appropriate, and identifying two nonzero vectors $v$, $w$ when
\begin{equation}
	w = A^j(v)
\end{equation}
for some integer $j$.  Thus $\mathcal{H}_A$ is a compact real or
complex manifold, according to whether one starts with ${\bf R}^n$ or
${\bf C}^n$.  There is a natural smooth mapping from ${\bf R}^n
\backslash \{0\}$ or ${\bf C}^n \backslash \{0\}$ onto
$\mathcal{H}_A$, as appropriate, in which a nonzero vector $v$ is
mapped to the corresponding equivalence class in $\mathcal{H}_A$, and
this mapping is holomorphic in the complex case.  Of course
$\mathcal{H}_A$ has the same dimension as ${\bf R}^n \backslash \{0\}$
or ${\bf C}^n \backslash \{0\}$, as appropriate, and the mapping to
$\mathcal{H}_A$ is a local diffeomorphism.

	If $n = 1$, then $A(v) = \alpha \, v$ for some nonzero scalar
$\alpha$ with $|\alpha| < 1$.  In the real case, if $\alpha > 0$, then
$A$ sends the positive real numbers to the positive real numbers and
$A$ sends the negative real numbers to the negative real numbers, and
$\mathcal{H}_A$ is basically a disjoint union of two circles.
If $\alpha < 0$, then $A$ maps the positive real numbers to the negative
real numbers and vice-versa, and $\mathcal{H}_A$ reduces to a single
circle.  In the complex case, we can think of ${\bf C} \backslash \{0\}$
as the same as ${\bf C} / 2 \pi i {\bf Z}$, i.e., the quotient of
${\bf C}$ by translations by integer multiples of $2 \pi i$, using the
exponential mapping from ${\bf C}$ onto ${\bf C} \backslash \{0\}$.
If $\beta$ is a complex number such that $\exp \beta = \alpha$ and
$L$ is the lattice in ${\bf C}$ consisting of complex numbers of the
form $2 m \pi i + n \beta$, $m, n \in {\bf Z}$, then $\mathcal{H}_A$
can be identified with ${\bf C} / L$, which is to say that we get a
$1$-dimensional complex torus.

	When $n \ge 2$, consider the special case where $A$ is of the
form $A(v) = \alpha \, v$ for some nonzero scalar $\alpha$ such that
$|\alpha| < 1$.  In this case there is a natural smooth mapping from
$\mathcal{H}_A$ onto ${\bf RP}^{n-1}$ or ${\bf CP}^{n-1}$, as
appropriate, in which elements of $\mathcal{H}_A$ are sent to the
lines in ${\bf R}^n$ or ${\bf C}^n$ that contain the corresponding
vectors.  This mapping is a submersion, and it is holomorphic in the
complex case.  The fibers of this mapping are copies of what one gets
in the $1$-dimensional case.  Also, linear mappings on ${\bf R}^n$
or ${\bf C}^n$, as appropriate, induce interesting mappings on 
$\mathcal{H}_A$.

	In general dimensions and for general $A$, suppose that $L$ is
a nontrivial linear subspace of ${\bf R}^n$ or ${\bf C}^n$, as
appropriate, such that
\begin{equation}
	A(L) = L.
\end{equation}
We can apply the same construction to get a space analogous to
$\mathcal{H}_A$ for the restriction of $A$ to $L$, and this space can
be viewed as a submanifold of $\mathcal{H}_A$.  In particular, a
$1$-dimensional invariant subspace for $A$ is the same as the span of
a nonzero eigenvector of $A$, and this leads to a submanifold of
$\mathcal{H}_A$ which is a copy of what one gets in the
$1$-dimensional case.  Also, a linear transformation on ${\bf R}^n$
or on ${\bf C}^n$ which commutes with $A$ leads to an interesting
mapping on $\mathcal{H}_A$.

	Suppose again that $M$, $N$ are smooth manifolds of dimensions
$m$, $n$, that $f : M \to N$ is a submersion, and that we have a connection
on $M$ associated to this submersion, defined by a smooth family $H_p$
of horizontal linear subspaces of $T_p M$, $p \in M$.  Fix a point
$p \in M$, and let $\alpha(t)$ be a smooth curve in $N$ defined on an
interval $[a, b]$ in the real line such that
\begin{equation}
	f(p) = \alpha(a).
\end{equation}
Consider the question of having a smooth curve $\beta(t)$, $a \le t \le b$,
in $N$ such that
\begin{equation}
	\beta(a) = p
\end{equation}
and
\begin{equation}
	f(\beta(t)) = \alpha(t) \quad\hbox{for all } a \le t \le b,
\end{equation}
so that $\beta(t)$ is a lifting of $\alpha(t)$ which starts at $p$,
and also
\begin{equation}
	{\dot \beta}(t) \in H_{\beta(t)}
				\quad\hbox{for all } a \le t \le b.
\end{equation}
Here ${\dot \beta}(t)$ denotes the derivative of $\beta(t)$ at $t$,
which is automatically an element of $T_{\beta(t)} M$, and which
satisfies
\begin{equation}
	df_{\beta(t)} ({\dot \beta}(t)) = {\dot \alpha}(t)
\end{equation}
for all $t$.  This condition and the requirement that the
derivative of $\beta(t)$ belong to the horizontal subspaces
specified by the connection determine the derivative of $\beta(t)$
in terms of $\beta(t)$ and the derivative of $\alpha(t)$.

	In other words, $\beta(t)$ satisfies an ordinary differential
equation.  More precisely, in local coordinates one can convert this
into a system of ordinary differential equations in ${\bf R}^m$ of the
usual type.  Standard results about ordinary differential equations
imply that $\beta(t)$ is uniquely determined by $\alpha(t)$ and
the starting point $p$, when such a lifting exists.  Also, such a
lifting always exists at least on a shorter interval beginning
at $a$.  If the submersion $f : M \to N$ is proper, then the
lifting $\beta(t)$ exists for the whole interval $[a, b]$, and
there are other conditions like this for lifting the whole curve
as well, basically by ensuring that the lifted curve remain in a
compact subset of $M$.

	Standard results about ordinary differential equations also
imply smoothness results about mappings associated to liftings like
these.  Namely, one can vary the choice of $p$ in the fiber
$f^{-1}(\alpha(a))$, and get smooth dependence on $p$ of the lifting.
One can always do this locally, for $t$ near $a$.  If $f : M \to N$ is
proper, so that we have liftings on the whole interval $[a, b]$, then
the lifting of paths defines a mapping from the fiber
$f^{-1}(\alpha(a))$ to the fiber $f^{-1}(\alpha(b))$, and this mapping
is a bijection with the inverse mapping obtained by running the
lifting backwards along the interval $[a, b]$.  Smooth dependence on
$p$ implies that this mapping from the fiber $f^{-1}(\alpha(a))$
onto $f^{-1}(\alpha(b))$ is in fact a diffeomorphism.

	If $f : M \to N$ is proper and $N$ is connected, then we get
another way to see that the fibers of $f$ are all diffeomorphic to
each other.  Namely, any pair of points in $N$ can be connected by a
smooth curve if $N$ is connected.  One can use liftings of this
curve to get a diffeomorphism between the corresponding fibers, as
above.  Of course there are plenty of variations of these themes.

	Let us consider another example.  Let $U$ denote the upper
half-space in ${\bf C}$, which consists of complex numbers with
positive imaginary part.  We start basically with the Cartesian 
product ${\bf C} \times U$, and the coordinate projection of
this space onto $U$.  This projection is obviously a holomorphic
mapping and a submersion.

	For each $\alpha \in U$, let $L_\alpha$ be the lattice in
${\bf C}$ consisting of $m + n \, \alpha$, $m, n \in {\bf Z}$.  We can
think of this as a fixed lattice in ${\bf C}$, or as a family of
lattices in a family of copies of ${\bf C}$.  For each $\alpha \in U$.
let us write $\mathcal{E}(\alpha)$ for the $1$-dimensional complex
torus ${\bf C} / L_\alpha$.  Let us write $\mathcal{E}$ for the space
that we get by identifying $(z, \alpha)$ and $(w, \alpha)$ in ${\bf C}
\times U$ when $w - z \in L_\alpha$.

	Thus $\mathcal{E}$ is a complex manifold of dimension $2$, and
we get a holomorphic projection mapping from ${\bf C} \times U$ onto
$\mathcal{E}$ sending a given pair $(z, \alpha)$ in ${\bf C} \times U$
to the corresponding equivalence class in $\mathcal{E}$.  Locally
$\mathcal{E}$ looks like ${\bf C} \times U$, which is to say that the
natural quotient mapping is locally a biholomorphism.  We also have a
natural mapping from $\mathcal{E}$ onto $U$, which is to say that the
standard coordinate projection from ${\bf C} \times U$ pushes down to
$\mathcal{E}$ in a natural way.  That is, the standard coordinate
projection from ${\bf C} \times U$ is the same as the natural
projection from ${\bf C} \times U$ onto $\mathcal{E}$ followed by the
mapping from $\mathcal{E}$ to $U$.  This mapping from $\mathcal{E}$ to
$U$ is a proper holomorphic submersion.

	For each $\alpha \in U$, we can identify the fiber in
$\mathcal{E}$ over $\alpha$ from our mapping $\mathcal{E} \to U$ with
$\mathcal{E}(\alpha)$ in a simple way.  As real manifolds, these
fibers are all diffeomorphic to each other.  However, these fibers are
not equivalent in general as complex manifolds.  In fact, for distinct
nearby $\alpha$'s the corresponding $\mathcal{E}(\alpha)$'s are not
holomorphically equivalent.  Locally $\mathcal{E}$ is holomorphically
equivalent to a product, and there is significant activity more
globally.

	Now let $M$, $N$ be smooth manifolds with dimensions
$m$, $n$, let $f : M \to N$ be a smooth submersion, and let
$H_p$, $p \in M$, be a smooth family of horizontal linear
subspaces of the tangent spaces of $M$ defining a connection
for the submersion.  Thus each $H_p$ has dimension $n$ and the
restriction of the differential $df_p$ of $f$ at $p$ to $H_p$
defines a one-to-one linear mapping onto $T_{f(p)} N$.  We
would like to describe the \emph{curvature} of this connection.
Let us first review some aspects of \emph{vector fields} on a
smooth manifold.

	Let $U$ be a nonempty open subset of $M$.  A smooth vector
field $X$ on $U$ assigns to each $p \in U$ a tangent vector
$X(p)$ to $M$ at $p$, in a way which is smooth in $p$.  Such a vector
field defines a first-order linear differential operator acting
on smooth real-valued functions on $U$, which is to say that if
$h$ is a smooth real-valued function on $U$, then $X(h)$ at a
point $p \in U$ is the directional derivative of $h$ in the direction
$X(p)$ at $p$.  These differential operators are linear, so that
\begin{equation}
	X(c_1 \, h_1 + c_2 \, h_2) 
		= c_1 \, X(h_1) + c_2 \, X(h_2)
\end{equation}
for all real numbers $c_1$, $c_2$ and all smooth functions $h_1$,
$h_2$, and they satisfy the Leibniz rule
\begin{equation}
	X(h_1 \, h_2) = X(h_1) \, h_2 + h_1 \, X(h_2),
\end{equation}
for differentiating the product of two smooth functions
$h_1$, $h_2$ on $U$.

	If $X_1$, $X_2$ are two smooth vector fields on $U$,
then one gets the associated Lie bracket $[X_1, X_2]$ of
$X_1$, $X_2$.  In terms of differential operators, we have
\begin{equation}
	[X_1, X_2] (h) = X_1(X_2(h)) - X_2(X_1(h))
\end{equation}
for all smooth functions $h$ on $U$.  Of course $X_1(X_2(h))$,
$X_2(X_1(h))$ involve second derivatives of $h$, and these second
derivatives cancel out in the difference, leaving a first-order
operator associated to a vector field.  If $X_1$, $X_2$ are smooth
vector fields on $U$ and $\phi_1$, $\phi_2$ are smooth real-valued
functions on $U$, then $\phi_1 \, X_1$, $\phi_2 \, X_2$ also define
smooth vector fields on $U$, and we have that
\begin{equation}
	[\phi_1 \, X_1, \phi_2 \, X_2]
		= \phi_1 \, \phi_2 \, [X_1, X_2]
			+ \phi_1 \, X_1(\phi_2) \, X_2
			- \phi_2 \, X_2(\phi_1) \, X_1.
\end{equation}

	Now let us return to the setting of our submersion and
connection, and assume that $X_1$, $X_2$ are smooth vector fields on a
nonempty open subset $U$ of $M$ such that $X_1(p)$, $X_2(p)$ are
elements of the horizontal linear subspace $H_p$ of the tangent space
$T_p M$ of $M$ at $p$ for each $p \in U$.  Define the curvature
$\mathcal{C}(X_1, X_2)$ at a point $p \in U$ to be the vertical
component of $[X_1, X_2]$ at $p$, which is to say that
$\mathcal{C}(X_1, X_2)$ at $p$ is an element of the vertical linear
subspace $V_p$ of the tangent space $T_p M$ of $M$ at $p$, and $[X_1,
X_2] - \mathcal{C}(X_1, X_2)$ lies in the horizontal subspace $H_p$
at $p$.  If $\phi_1$, $\phi_2$ are smooth real-valued functions
on $U$, then $\phi_1 \, X_1$, $\phi_2 \, X_2$ are also horizontal
vector fields on $U$, and
\begin{equation}
	\mathcal{C}(\phi_1 \, X_1, \phi_2 \, X_2)
		= \phi_1 \, \phi_2 \, \mathcal{C}(X_1, X_2).
\end{equation}
This shows that the curvature $\mathcal{C}(X_1, X_2)$, at a point $p
\in U$, depends only on the values of $X_1$, $X_2$ at $p$, and thus at
$p$ the curvature $\mathcal{C}(X_1, X_2)$ defines an antisymmetric
bilinear mapping from $H_p \times H_p$ into $V_p$, which depends
smoothly on $p$.  Using the isomorphism between $H_p$ and $T_{f(p)} N$
given by $df_p$, one can reformulate the curvature by saying that for
each $p \in M$ it is an antisymmetric bilinear mapping from $T_{f(p)}
N \times T_{f(p)} N$ into $V_p$ which depends smoothly on $p$.

	What does it mean for the curvature to be equal to $0$
everywhere on $M$?  This is equivalent to saying that if $X_1$, $X_2$
are smooth horizontal vector fields on a nonempty subset $U$ of $M$,
then the Lie bracket $[X_1, X_2]$ is also a horizontal vector field on
$U$.  In other words, the curvature of the connection is equal to $0$
on all of $M$ if and only if the corresponding distribution of
horizontal linear subspaces of the tangent spaces is integrable.  By a
well-known theorem of Frobenius, this means that there is a foliation
of $M$ by $n$-dimensional smooth submanifolds whose tangent spaces are
exactly the horizontal linear subspaces of the tangent spaces of $M$
given by the connection.

\section{Metric spaces}
\label{metric spaces}
\setcounter{equation}{0}

	By a \emph{metric space} we mean a nonempty set $M$ together
with a real-valued function $d(x, y)$ defined for $x, y \in M$, called
the \emph{distance function} or \emph{metric} on $M$, such that $d(x,
y) \ge 0$ for all $x, y \in M$, $d(x, y) = 0$ if and only if $x = y$,
\begin{equation}
	d(y, x) = d(x, y)
\end{equation}
for all $x, y \in M$, and
\begin{equation}
	d(x, z) \le d(x, y) + d(y, z)
\end{equation}
for all $x, y, z \in M$.  A basic example if given by the real numbers
${\bf R}$ equipped with the standard metric $|x - y|$.  Recall that if
$x$ is a real number, then the \emph{absolute value} of $x$ is denoted
$|x|$ and defined to be equal to $x$ when $x \ge 0$ and to $-x$ when
$x \le 0$, and that
\begin{equation}
	|x + y| \le |x| + |y|, \quad |x \, y| = |x| \, |y|
\end{equation}
for all $x, y \in {\bf R}$.

	If $(M, d(x,y))$ and $(N, \rho(u, v))$ are metric spaces
and $f : M \to N$ is a mapping from $M$ to $N$, then $f$ is said
to be \emph{continuous} if for every $x \in M$ and every positive
real number $\epsilon$ there is a positive real number $\delta$
such that
\begin{equation}
	\rho(f(y), f(x)) < \epsilon
\end{equation}
for all $y \in M$ such that $d(y, x) < \delta$.  For instance,
constant mappings are always continuous, and the identity mapping on a
metric space $(M, d(x,y))$ is continuous as a mapping from $M$ to
itself.  More generally, if $(M, d(x,y))$ is a metric space and $E$ is
a nonempty subset of $M$, then we can consider $E$ to be a metric
space itself, using the same metric $d(x, y)$ restricted to $E$, and
then the inclusion mapping of $E$ into $M$, which sends each element
of $E$ to itself, is continuous as a mapping from $E$ to $M$.

	Let $(M, d(x, y))$ be a metric space, and let $p$ be a point
in $M$.  From the triangle inequality it is easy to see that
\begin{equation}
	d(x, p) - d(y, p) \le d(x, y)
\end{equation}
for all $x, y \in M$, and similarly
\begin{equation}
	d(y, p) - d(x, p) \le d(x, y),
\end{equation}
so that
\begin{equation}
	|d(x, p) - d(y, p)| \le d(x, y)
\end{equation}
for all $x, y \in M$.  It follows that the real-valued function
$f_p(x) = d(x, p)$ on $M$ is continuous, where, as usual, we employ
the standard metric on the real numbers.

	If $f_1$, $f_2$ are two real-valued continuous functions on a
metric space $(M, d(x, y))$, then the sum $f_1 + f_2$ and the product
$f_1 \, f_2$ are also continuous functions.  This is not too difficult
to check.  Similarly, if $f$ is a continuous real-valued function on
$M$ such that $f(x) \ne 0$ for all $x \in M$, then $1/f(x)$ is also a
continuous function on $M$.

	Let $(M, d(x,y))$ be a metric space, and let $A$ be a nonempty
subset of $M$.  We denote the distance from a point $x$ in $M$ to $A$
by $\dist(x, A)$, and we define it to be the infimum of $d(x, y)$ over
all $y \in A$.  If $A$, $B$ are two nonempty subsets of $M$, then the
distance between $A$ and $B$ is denoted $\dist(A, B)$ and defined to
be the infimum of $d(x, y)$ over all $x \in A$ and $y \in B$.

	If $A$ is a nonempty subset of $M$ and $x$, $y$ are elements
of $M$, then one can check that
\begin{equation}
	\dist(x, A) \le \dist(y, A) + d(x, y).
\end{equation}
As a result,
\begin{equation}
	|\dist(x, A) - \dist(y, A)| \le d(x, y)
\end{equation}
for all $x, y \in M$.  In particular, $\dist(x, A)$ is a real-valued
continuous function of $x$ on $M$.

	If $A$ and $B$ are nonempty subsets of $M$ and $t$ is a
positive real number, then we say that $A$, $B$ are $t$-close if for
each $a \in A$ there is a $b \in B$ such that $d(a, b) < t$, and if
for each $b \in B$ there is an $a \in A$ such that $d(a, b) < t$.  If
$A$, $B$ are nonempty subsets of $M$ which are $t$-close for some
positive real number $t$, then the Hausdorff distance from $A$ to $B$
is denoted $D(A, B)$ and defined to be the infimum of the positive
real numbers $t$ such that $A$, $B$ are $t$-close.  In this case, if
$x$ is any point in $M$, then
\begin{equation}
	\dist(x, A) \le \dist(x, B) + D(A, B).
\end{equation}

	A subset $A$ of $M$ is said to be \emph{bounded} if there is a
point $p$ in $M$ and a positive real number $R$ such that $d(x, p) \le
R$ for all $x \in A$.  It is easy to see that once this holds for some
$p \in M$, it works for all $p \in M$, with a choice of $R$ that
depends on $p$.  The \emph{diameter} of a nonempty bounded subset $A$
of $M$ is denoted $\diam A$ and defined to be the supremum of $d(x,
y)$ over all $x, y \in A$.

	If $A$, $B$ are bounded subsets of $M$, then clearly $A$, $B$
are $t$-close for some positive real numbers $t$, and thus the
Hausdorff distance $D(A, B)$ is defined.  Of course $B$, $A$ are
$t$-close when $A$, $B$ are $t$-close, so that $D(A, B) = D(B, A)$.
Also, if $A$, $B$, $C$ are nonempty subsets of $M$ and $s$, $t$ are
positive real numbers such that $A$, $B$ are $s$-close and $B$, $C$
are $t$-close, then $A$, $C$ are $s + t$ close, and 
\begin{equation}
	D(A, C) \le D(A, B) + D(B, C).
\end{equation}

	Suppose that $\phi(x)$ is a monotone increasing real-valued
function on the real line, so that $\phi(x) \le \phi(y)$ when $x$, $y$
are real numbers such that $x \le y$.  For each real number $x$, the
left and right-sided limits of $\phi$ at $x$, denoted $\phi(x-)$ and
$\phi(x+)$, respectively, automatically exist and can be given by
\begin{equation}
	\phi(x-) = \sup \{\phi(w) : w < x\}, \quad
		\phi(x+) = \inf \{\phi(y) : y > x\}.
\end{equation}
Clearly
\begin{equation}
	\phi(x-) \le \phi(x) \le \phi(x+),
\end{equation}
and $\phi$ is continuous at $x$ if and only if
\begin{equation}
	\phi(x-) = \phi(x+).
\end{equation}

	For each positive real number $p$, one can show that the
function $|x|^p$ is a continuous real-valued function on ${\bf R}$.
Fix a positive integer $n$, and for each positive real number $p$
consider the real-valued function $\|x\|_p$ on ${\bf R}^n$ defined by
\begin{equation}
	\|x\|_p = \biggl(\sum_{j=1}^n |x_j|^p \biggr)^{1/p},
\end{equation}
$x = (x_1, \ldots, x_n)$.  We can also allow $p = \infty$ here
by setting
\begin{equation}
	\|x\|_\infty = \max (|x_1|, \ldots, |x_n|).
\end{equation}

	Thus $\|x\|_p$ is a nonnegative real number for each $x \in
{\bf R}^n$ and $0 < p \le \infty$ which is equal to $0$ if and only if
$x = 0$.  We also have that
\begin{equation}
	\|t \, x\|_p = |t| \, \|x\|_p
\end{equation}
for each real number $t$, $x \in {\bf R}^n$, and $0 < p \le \infty$.
Here $t \, x$ denotes the usual scalar multiplication of $t$ and
the vector $x$, so that
\begin{equation}
	t \, x = (t \, x_1, \ldots, t \, x_n).
\end{equation}

	Clearly
\begin{equation}
	\|x\|_\infty \le \|x\|_p
\end{equation}
for all $x \in {\bf R}^n$ and $0 < p < \infty$.  More generally, if
$0 < p \le q < \infty$, then
\begin{equation}
	\|x\|_q \le \|x\|_p.
\end{equation}
Indeed,
\begin{eqnarray}
	\|x\|_q^q = \sum_{j=1}^n |x_j|^q
		& \le & \|x\|_\infty^{q-p} \sum_{j=1}^n |x_j|^p		\\
		& = & \|x\|_\infty^{q-p} \, \|x\|_p^p
		\le \|x\|_p^q.					\nonumber
\end{eqnarray}

	When $0 < p \le 1$ we can take $q = 1$ and $n = 2$ to obtain
that
\begin{equation}
	(a + b)^p \le a^p + b^p
\end{equation}
for all nonnegative real numbers $a$, $b$.  For $p \ge 1$ a natural
counterpart of this is the fact that $t^p$ is a convex function of
$t$ on the set of nonnegative real numbers.  In other words,
if $t$, $u$ are nonnegative real numbers and $\lambda$ is a real
number such that $0 < \lambda < 1$, then
\begin{equation}
	(\lambda \, t + (1 - \lambda) \, u)^p
		\le \lambda \, t^p + (1 - \lambda) \, u^p
\end{equation}
for every real number $p \ge 1$.

	It is easy to see that
\begin{equation}
	\|x\|_p \le n^{1/p} \, \|x\|_\infty
\end{equation}
for every $x \in {\bf R}^n$ and every positive real number $p$.
In fact, if $p$, $q$ are positive real numbers such that $p \le q$,
then
\begin{equation}
	\|x\|_p \le n^{(1/p) - (1/q)} \, \|x\|_q
\end{equation}
for all $x \in {\bf R}^n$.  This can be derived from the convexity of
the function $t^{q/p}$ for $t \ge 0$.

	For any $x, y \in {\bf R}^n$ and $1 \le p \le \infty$ we have
that
\begin{equation}
	\|x + y\|_p \le \|x\|_p + \|y\|_p.
\end{equation}
This is easy to derive directly from the definitions when $p = 1$ or
$\infty$.  In general, one can use homogeneity to reduce to showing
that
\begin{equation}
	\|\lambda \, x + (1 - \lambda) \, y\|_p \le 1
\end{equation}
when $x, y \in {\bf R}^n$ satisfy $\|x\|_p, \|y\|_p \le 1$ and
$\lambda$ is a real number such that $0 < \lambda < 1$, and for $p \ge
1$ this can be derived from the convexity of $t^p$, $t \ge 0$.

	As a result, when $1 \le p \le \infty$, we have that
\begin{equation}
	d_p(x, y) = \|x - y\|_p
\end{equation}
defines a metric on ${\bf R}^n$.  When $0 < p \le 1$ we can set
\begin{equation}
	d_p(x, y) = \|x - y\|_p^p,
\end{equation}
and this also defines a metric on ${\bf R}^n$.  To be more precise,
this uses the fact that
\begin{equation}
	\|v + w\|_p^p \le \|v\|_p^p + \|w\|_p^p
\end{equation}
when $v, w \in {\bf R}^n$ and $0 < p \le 1$.

	If $(M, d(x,y))$ and $(N, \rho(u, v))$ are metric spaces
and $f : M \to N$ is a mapping from $M$ into $N$, then we say that
$f$ is \emph{Lipschitz} if there is a nonnegative real number $C$
such that
\begin{equation}
	\rho(f(x), f(y)) \le C \, d(x, y)
\end{equation}
for all $x, y \in M$.  We might say that $f$ is $C$-Lipschitz in this
case, and notice that $0$-Lipschitz mappings are constant.  Of course
Lipschitz mappings are automatically continuous.

	Suppose that $(M_1, d_1(x, y))$, $(M_2, d_2(u, v))$, and
$(M_3, d_3(z, w))$ are metric spaces, and that $f_1 : M_1 \to M_2$ and
$f_2 : M_2 \to M_3$ are mappings between them.  As usual, the
composition $f_2 \circ f_1$ is the mapping from $M_1$ to $M_3$ defined
by
\begin{equation}
	(f_2 \circ f_1)(x) = f_2(f_1(x))
\end{equation}
for all $x \in M$.  It is easy to see that if $f_1$ is a continuous
mapping from $M_1$ to $M_2$ and $f_2$ is a continuous mapping from
$M_2$ to $M_3$, then the composition $f_2 \circ f_1$ is a continuous
mapping from $M_1$ to $M_3$, and that if $f_1$, $f_2$ are Lipschitz,
then so is the composition $f_2 \circ f_1$.

	It is easy to generate examples of real-valued Lipschitz
functions on the real line, which can then be composed with some of
the basic real-valued Lipschitz functions on a metric space mentioned
earlier to produce more Lipschitz functions.  In general, if $f_1$,
$f_2$ are two real-valued Lipschitz functions on a metric space $M$,
then $f_1 + f_2$, $\min(f_1, f_2)$, and $\max(f_1, f_2)$ are also
Lipschitz functions on $M$, and a real number times a real-valued
Lipschitz function is again a Lipschitz function.  For products of
Lipschitz functions, or reciprocals of nonzero Lipschitz functions,
the situation is more complicated, although there are simple
sufficient conditions for the result to be Lipschitz.

	Now let us look at continuous curves or paths in a metric
space.  Namely, let $(M, d(x,y))$ be a metric space, and let $I$
be a closed interval in the real line.  That is, $I$ might be of
the form 
\begin{equation}
	[a, b] = \{x \in {\bf R} : a \le x \le b\}
\end{equation}
for some real numbers $a$, $b$ with $a \le b$, in which case
$I$ is a closed and bounded interval, or $I$ might be an
unbounded closed interval, of the form
\begin{equation}
	[a, \infty) = \{x \in {\bf R} : x \ge a\}
\end{equation}
for some real number $a$, or
\begin{equation}
	(-\infty, b] = \{x \in {\bf R} : x \le b\}
\end{equation}
for some real number $b$, or
\begin{equation}
	(-\infty, \infty) = {\bf R}.
\end{equation}

	A continuous path in $M$ parameterized by the interval
$I$ is simply a continuous mapping from $I$ into $M$.  Sometimes
we are particularly interested in paths which are Lipschitz.
This can be interpreted as meaning that the path has bounded
speed.

	Suppose that $I = [a, b]$ is a closed and bounded interval
in the real line, and that $p : I \to M$ is a continuous path 
on $I$ in the metric space $(M, d(x, y))$.  By a partition of
$I$ we mean a finite sequence $\mathcal{P} = \{t_j\}_{j=0}^k$
of real numbers such that
\begin{equation}
	a = t_0 < \ldots < t_k = b.
\end{equation}
Associated to this partition $\mathcal{P}$ we get an approximation
to the length of the path $p$, defined by
\begin{equation}
	\Lambda_a^b(p, \mathcal{P}) = \sum_{j=1}^k d(p(t_j), p(t_{j-1})).
\end{equation}

	If this quantity is uniformly bounded over all partitions
$\mathcal{P}$ of $I$, then we say that the path $p$ has finite length,
and we define the length of $p$, denoted $\Lambda_a^b(p)$, to be the
supremum of $\Lambda_a^b(p, \mathcal{P})$ over all partitions
$\mathcal{P}$ of $I$.  If $p : I \to M$ is $C$-Lipschitz for some
nonnegative real number $C$, then $p$ has finite length, and
\begin{equation}
	\Lambda_a^b(p) \le C \, (b - a).
\end{equation}
It is sometimes convenient to allow $a = b$ and $k = 0$ in the
definition of a partition, in which case the path automatically
has length $0$, and in general a path has length $0$ if and only
if it is constant.

	If $(M, d(x, y))$ is a metric space, $p$ is an element of $M$,
and $r$ is a positive real number, then the \emph{open ball in $M$
with center $p$ and radius $r$} is denoted $B(p, r)$ and defined by
\begin{equation}
	B(p, r) = \{z \in M : d(p, z) < r\}.
\end{equation}
Similarly, the \emph{closed ball with center $p$ and radius $r$}
is denoted $\overline{B}(p, r)$ and defined by
\begin{equation}
	\overline{B}(p, r) = \{z \in M : d(p, z) \le r\}.
\end{equation}
We make the convention that a ``ball'' means an open ball unless
otherwise specified.

	A subset $U$ of $M$ is said to be \emph{open} if for every
point $p \in U$ there is a positive real number $r$ such that
\begin{equation}
	B(p, r) \subseteq U.
\end{equation}
The union of any family of open sets is open, and the intersection of
finitely many open sets is open.  Note that the empty set $\emptyset$
and $M$ itself are automatically open subsets of $M$, and one can
check that open balls in $M$ are open subsets of $M$.

	A sequence $\{x_j\}_{j=1}^\infty$ of points in $M$ is said to
\emph{converge} to a point $x$ in $M$ if for every $\epsilon > 0$
there is a positive integer $L$ such that
\begin{equation}
	d(x_j, x) < \epsilon
\end{equation}
for all $j \ge L$.  In this case we write
\begin{equation}
	\lim_{j \to \infty} x_j = x,
\end{equation}
and we call $x$ the limit of the sequence $\{x_j\}_{j=1}^\infty$.
It is not difficult to see that the limit of a sequence is unique
if it exists.

	A subset $F$ of $M$ is said to be \emph{closed} if every
sequence of points in $F$ which converges to some point in $M$
has its limit in $F$.  The empty set and $M$ itself are automatically
closed sets, and one can check that closed balls in $M$ are closed
subsets of $M$.  The intersection of any family of closed sets is
closed, and the union of finitely many closed sets is again closed.

	In fact, a subset $U$ of $M$ is open if and only if its
complement $M \backslash U$ is closed.  Recall that the complement
of a subset $E$ of $M$ in $M$ is given by
\begin{equation}
	M \backslash E = \{x \in M : x \not\in E\}.
\end{equation}
Equivalently, a subset $F$ of $M$ is closed if and only if $M
\backslash F$ is an open subset of $M$.

	Suppose that $M$, $N$ are sets and that $f$ is a mapping from
$M$ to $N$.  If $A$ is a subset of $M$, then the image of $A$ under
$f$ is denoted $f(A)$ and is the subset of $N$ defined by
\begin{equation}
	f(A) = \{f(x) : x \in A\}.
\end{equation}
In particular, the image of $f$ simply means $f(M)$.  If $B$ is a
subset of $N$, then the inverse image of $B$ under $f$ is denoted
$f^{-1}(B)$ and is the subset of $M$ defined by
\begin{equation}
	f^{-1}(B) = \{x \in M : f(x) \in B\}.
\end{equation}

	The image of the union of a family of subsets of $M$ under $f$
is equal to the union of the images of the individual subsets, and the
image of the intersection of a family of subsets of $M$ under $f$
is contained in the intersection of the images of the individual subsets.
The inverse image of the union of a family of subsets of $N$ under
$f$ is equal to the union of the inverse images of the individual
subsets of $N$, and the inverse image of the intersection of a family
of subsets of $N$ is equal to the intersection of the corresponding
individual inverse images.  If $B$ is a subset of $N$, then
\begin{equation}
	M \backslash f^{-1}(B) = f^{-1}(N \backslash B).
\end{equation}

	If $(M, d(x, y))$ and $(N, \rho(u, v))$ are metric spaces,
then a mapping $f$ from $M$ to $N$ is continuous if and only if
$f^{-1}(V)$ is an open subset of $M$ for every open subset $V$ of $N$.
This is also equivalent to saying that $f^{-1}(E)$ is a closed subset
of $M$ for every closed subset $E$ of $N$.  Moreover, $f$ is
continuous if and only if for every sequence $\{x_j\}_{j=1}^\infty$ of
points in $M$ which converges to a point $x$ in $M$ we have that the
sequence $\{f(x_j)\}_{j=1}^\infty$ converges to $f(x)$ in $N$.

	A sequence $\{x_j\}_{j=1}^\infty$ of points in a metric space
$(M, d(y, z))$ is said to be a \emph{Cauchy sequence} if for every
$\epsilon > 0$ there is a positive integer $L$ such that
\begin{equation}
	d(x_j, x_k) < \epsilon
\end{equation}
for all $j, k \ge L$.  Every convergent sequence is a Cauchy sequence,
and a metric space is said to be \emph{complete} if every Cauchy
sequence in the space converges to some point in the space.
A basic property of Euclidean spaces ${\bf R}^n$ with their
standard metrics $\|x - y\|_2$ is that they are complete.

	If $(M, d(x, y))$ and $(N, \rho(u, v))$ are metric spaces and
$f$ is a mapping from $M$ to $N$, then $f$ is said to be
\emph{uniformly continuous} if for every $\epsilon > 0$ there is a
$\delta > 0$ such that
\begin{equation}
	\rho(f(x), f(y)) < \epsilon
\end{equation}
for all $x, y \in M$ such that $d(x, y) < \delta$.  It is easy to
see that Lipschitz mappings are uniformly continuous.  Constant
mappings are uniformly continuous trivially, and the identity mapping
on a metric space is $1$-Lipschitz and hence uniformly continuous.

	If $f : M \to N$ is uniformly continuous and
$\{x_j\}_{j=1}^\infty$ is a Cauchy sequence in $M$, then
$\{f(x_j)\}_{j=1}^\infty$ is a Cauchy sequence in $N$.  In particular,
if $N$ is complete, then $\{f(x_j)\}_{j=1}^\infty$ converges in $N$.
Also, if $f : M \to N$ is uniformly continuous and
$\{x_j\}_{j=1}^\infty$, $\{y_j\}_{j=1}^\infty$ are sequences in $M$
such that
\begin{equation}
	\lim_{j \to \infty} d(x_j, y_j) = 0,
\end{equation}
then
\begin{equation}
	\lim_{j \to \infty} \rho(f(x_j), f(y_j)) = 0
\end{equation}
too.

	In any metric space $(M, d(x, y))$, the \emph{closure} of a
subset $E$ is denoted $\overline{E}$ can be defined as the set of
points $x \in M$ for which there is a sequence $\{x_j\}_{j=1}^\infty$
of points in $E$ which converges to $x$.  Thus
\begin{equation}
	E \subseteq \overline{E}
\end{equation}
automatically, and one can check that $\overline{E}$ is a closed
subset of $M$ which is contained in any other closed subset of $M$
that contains $E$.  A subset $E$ of $M$ is said to be \emph{dense}
in $M$ if $\overline{E} = M$.

	Suppose that $(M, d(x, y))$ and $(N, \rho(u, v))$ are
metric spaces, with $N$ complete, $E$ is a dense subset of $M$,
and $f$ is a uniformly continuous mapping from $E$ to $N$.
Under these conditions, one can show that there is a uniformly
continuous mapping from $M$ to $N$ which agrees with $f$ on $E$.
This extension is unique, and for that matter if $f_1$, $f_2$
are two continuous mappings from one metric space into another,
then the set of points in the domain where $f_1$ and $f_2$ agree
is a closed set.

	If $\{x_j\}_{j=1}^\infty$ is a sequence of points in some set
$A$, and if $\{j_k\}_{k=1}^\infty$ is a strictly increasing sequence
of positive integers, then $\{x_{j_k}\}_{k=1}^\infty$ is called a
\emph{subsequence} of the original sequence $\{x_j\}_{j=1}^\infty$.  A
subset $K$ of a metric space $(M, d(x, y))$ is said to be
\emph{compact} if every sequence of points in $K$ has a subsequence
which converges to a point in $K$.  Notice that if $Y$ is a nonempty
subset of $M$ such that $K \subset Y$, then $K$ is compact as a subset
of $M$ if and only if $K$ is compact as a subset of $Y$, viewed as a
metric space on its own, using the restriction of the metric from $M$.

	A compact subset of a metric space is always closed, basically
because any subsequence of a convergent sequence converges to the same
limit.  Notice that a Cauchy sequence in a metric space which has a
convergent subsequence also converges to the same limit.  As a result,
a Cauchy sequence contained in a compact subset of a metric space
converges to a point in that subset.

	A compact subset $K$ of a metric space $(M, d(x,y))$ is
bounded.  To see this, let $p$ be a point in $M$, and assume for the
sake of a contradiction that $K$ is not bounded, so that for each
positive integer $j$ there is a point $x_j \in K$ with $d(p, x_j) \ge
j$.  It is easy to see that the sequence $\{x_j\}_{j=1}^\infty$ cannot
have a convergent subsequence in this case, contradicting the
assumption that $K$ is compact.

	A subset $E$ of a metric space $(M, d(x, y))$ is said to be
\emph{totally bounded} if for every positive real number $r$
there is a finite set $F \subseteq E$ such that
\begin{equation}
	E \subseteq \bigcup_{x \in F} B(x, r).
\end{equation}
A compact subset of $M$ is also totally bounded.  Indeed, a subset
$E$ of $M$ is not totally bounded if and only if there is an
$\epsilon > 0$ and a sequence of points $\{x_j\}_{j=1}^\infty$
in $E$ such that $d(x_j, x_k) \ge \epsilon$ for all positive
integers $j$, $k$ with $j \ne k$.

	More precisely, a subset $E$ of a metric space $(M, d(x,y))$
is totally bounded if and only if every sequence of points in $E$ has
a subsequence which is a Cauchy sequence.  This is not too difficult
to show.  As a result, a subset of a complete metric space is compact
if and only if it is closed and totally bounded.  In particular,
closed and bounded subsets of Euclidean spaces are compact.

	Let $(M, d(x, y))$ and $(N, \rho(u, v))$ be metric spaces, and
let $f$ be a mapping from $M$ to $N$.  We say that $f$ is
\emph{bounded} if the image of $f$ is a bounded subset of $N$.  The
space of bounded continuous mappings from $M$ to $N$ is denoted
$C_b(M, N)$.

	There is a natural metric on $C_b(M, N)$, called the
\emph{supremum metric}, which is defined by
\begin{equation}
	\sigma(f_1, f_2) = \sup \{\rho(f_1(x), f_2(x)) : x \in M\}
\end{equation}
for $f_1, f_2 \in C_b(M, N)$.  Convergence of sequences in $C_b(M, N)$
is equivalent to \emph{uniform convergence}, as compared to pointwise
convergence of mappings.  A basic result, which is not too difficult
to show, states that if $N$ is a complete metric space, then so is
$C_b(M, N)$.

	Let us write $\Lip(M, N)$ for the space of Lipschitz mappings
from $M$ to $N$, and for each positive real number $k$, let us write
$\Lip_k(M, N)$ for the space of $k$-Lipschitz mappings from $M$ to
$N$.  If $M$ is bounded, then $\Lip(M, N)$ is contained in $C_b(M,
N)$, and for each $k > 0$, $\Lip_k(M, N)$ is a closed subset of
$C_b(M, N)$.  If $M$ is bounded, $p$ is an element of $M$, $k$
is a positive real number, and $B$ is a bounded subset of $N$,
then the set of $f$ in $\Lip_k(M, N)$ such that $f(p) \in B$
is a bounded subset of $C_b(M, N)$.

	Suppose that $(M, d(x, y))$ and $(N, \rho(u, v))$ are metric
spaces, with $M$ compact.  If $f$ is a continuous mapping from $M$ to
$N$, then the image of $f$ is a compact subset of $N$.  In particular,
every continuous mapping from $M$ to $N$ is bounded in this case.  If
$M$, $N$ are both compact, then one can show that $\Lip_k(M, N)$ is a
compact subset of $C_b(M, N)$ for every positive real number $k$.

\end{document}